\documentclass{article}
\usepackage[linesnumbered,ruled,vlined]{algorithm2e}
\usepackage[margin=1in]{geometry}
\usepackage[table]{xcolor}
\usepackage{amsmath,amsfonts,amssymb,url,array,textcomp,stfloats,url,verbatim,graphicx,cite,amsthm,soul,tikz,float,subcaption,multirow,doi,appendix,pdfpages}
\usepackage[textsize=small,textwidth=0.75in]{todonotes}
\usepgflibrary{shapes}

\newcommand\Renyi[0]{R{\'e}nyi }
\newcommand\ER[0]{Erd\H{o}s-\Renyi}
\newcommand\pdn[1]{\gamma_P\left(#1\right)}

\newcommand{\observed}[2]{\operatorname{Obs}\left(#1;#2\right)}
\newcommand{\supp}[1]{\partial\left(#1\right)}
\newcommand{\pref}[1]{\operatorname{Pref}\left(#1\right)}
\newcommand{\jeff}[1]{\operatorname{Pev}\left(#1\right)}

\definecolor{johns}{HTML}{c875c4}

\definecolor{beths}{HTML}{acfffc}

\newtheorem{theorem}{Theorem}[section]

\newtheorem{proposition}[theorem]{Proposition}
\newtheorem{observation}[theorem]{Observation}

\SetKwInput{KwInput}{Input}
\SetKwInput{KwOutput}{Output}

\tikzstyle{Black Outline Node}=[fill=white, draw=black, shape=circle, minimum size=1.5em]
\tikzstyle{Small Black Outline Node}=[fill=white, draw=black, shape=circle]
\tikzstyle{Meta Rect}=[fill=none, draw=none, shape=rectangle, inner sep=0pt, outer sep=0pt, text height=0pt]
\tikzstyle{Yellow Node}=[fill=yellow, draw=black, shape=circle]
\tikzstyle{Blue Node}=[fill=blue, draw=black, shape=circle]
\tikzstyle{Red Node}=[fill=red, draw=black, shape=circle]
\tikzstyle{Text Node}=[fill=none, draw=none, shape=rectangle, inner sep=0pt, outer sep=0pt, text height=9pt]
\tikzstyle{empty circle}=[fill=white, draw=black, shape=circle]
\tikzstyle{empty triangle}=[fill=white, draw=black, regular polygon, regular polygon sides=3]
\tikzstyle{partial circle}=[fill=black!25, draw=black, shape=circle]
\tikzstyle{solid circle}=[fill=black, draw=black, shape=circle]
\tikzstyle{solid star}=[fill=black, draw=black, shape=star]
\tikzstyle{partial diamond}=[fill=black!50, draw=black, shape=diamond]
\tikzstyle{Black Edge}=[-]
\tikzstyle{Directed Edge}=[->]
\tikzstyle{Red Edge}=[-, draw=red]
\tikzstyle{solid edge}=[-]
\tikzstyle{dashed edge}=[dashed]

\begin{document}

\title{The Power Domination Toolbox\thanks{This project was sponsored by the Air Force Research Laboratory via the Autonomy Technology Research Center}}
\date{\today}
\author{Johnathan Koch\thanks{Youngstown State University, Applied Research Solutions}\;, Beth Bjorkman\thanks{Air Force Research Laboratory}}

\maketitle
\begin{figure}[H]
    \centering
    \resizebox{!}{8cm}{
        \begin{tikzpicture}
            \node [style=Meta Rect] (8) at (0.75, 1) {};
            \node [style=Meta Rect] (13) at (3.25, 1.75) {};
            \node [style=Meta Rect] (14) at (3.25, 0.25) {};
            \node [style=Meta Rect] (26) at (3.25, 2.75) {};
            \node [style=Meta Rect] (27) at (0.75, 2.75) {};
            \node [style=Meta Rect] (35) at (1.5, 4.25) {};
            \node [style=Meta Rect] (41) at (2, 5.25) {};
            \node [style=Meta Rect] (42) at (2, 6.25) {};
            \node [style=Meta Rect] (43) at (2.5, 5.75) {};
            \node [style=Meta Rect] (44) at (1.5, 5.75) {};
            \node [style=Meta Rect] (82) at (0.5, 8) {};
            \node [style=Meta Rect] (83) at (0.5, 9) {};
            \node [style=Meta Rect] (88) at (2, 9) {};
            \node [style=Meta Rect] (89) at (2, 8) {};
            \node [style=Meta Rect] (90) at (2.5, 6.75) {};
            \node [style=Meta Rect] (91) at (2.5, 7.75) {};
            \node [style=Meta Rect] (96) at (3.5, 7.75) {};
            \node [style=Meta Rect] (97) at (3.5, 6.75) {};
            \node [style=Meta Rect] (98) at (1.25, 8) {};
            \node [style=Meta Rect] (99) at (2, 7.25) {};
            \node [style=Small Black Outline Node] (10) at (2, 1.75) {};
            \node [style=Yellow Node] (15) at (2.5, 1) {};
            \node [style=Yellow Node] (16) at (3.25, 1) {};
            \node [style=Blue Node] (17) at (1.5, 1) {};
            \node [style=Yellow Node] (22) at (1.5, 2.75) {};
            \node [style=Yellow Node] (23) at (2.5, 2.75) {};
            \node [style=Small Black Outline Node] (24) at (2, 3.25) {};
            \node [style=Small Black Outline Node] (25) at (2, 2.25) {};
            \node [style=Yellow Node] (32) at (2, 4.25) {};
            \node [style=Small Black Outline Node] (33) at (2.5, 4.75) {};
            \node [style=Small Black Outline Node] (34) at (2.5, 3.75) {};
            \node [style=Yellow Node] (40) at (2, 5.75) {};
            \node [style=Red Node] (84) at (1, 8.5) {};
            \node [style=Red Node] (85) at (1.5, 8.5) {};
            \node [style=Small Black Outline Node] (86) at (0.5, 8.5) {};
            \node [style=Small Black Outline Node] (87) at (2, 8.5) {};
            \node [style=Red Node] (92) at (3, 7.25) {};
            \node [style=Small Black Outline Node] (94) at (2.5, 7.25) {};
            \node [style=Small Black Outline Node] (95) at (3.5, 7.25) {};
            \node [style=Meta Rect] (100) at (-0.75, 11.25) {};
            \node [style=Text Node, text width = 3cm] at (-3,5.5) {Minimum Power Dominating Set Solution Space};
            \node [style=Text Node, label=right:{Graph Contraction}] (113) at (4.5, 8) {};
            \node [style=Text Node, label=right:{High Degree Vertices}] (114) at (4.5, 5.75) {};
            \node [style=Text Node, label=right:{Terminal Forts}] (115) at (4.5, 4.25) {};
            \node [style=Text Node, label=right:{Non-terminal Forts}] (116) at (4.5, 2.75) {};
            \node [style=Text Node, label=right:{Only Active Vertices}] (117) at (4.5, 1) {};
            \draw [style=Black Edge, bend left=40] (3.25, 2) to (4, 1);
            \draw [style=Black Edge, bend left=40] (4, 1) to (3.25, 0);
            \draw [style=Black Edge, bend right=40] (0.75, 2) to (0, 1);
            \draw [style=Black Edge, bend right=40] (0, 1) to (0.75, 0);
            \draw [style=Black Edge, bend right=40] (0.75, 3.5) to (0, 2.25);
            \draw [style=Black Edge, bend left=40] (3.25, 3.5) to (4, 2.25);
            \draw [style=Black Edge, bend right=40] (0.75, 5) to (0, 3.75);
            \draw [style=Black Edge, bend left=40] (3.25, 5) to (4, 3.75);
            \draw [style=Black Edge, bend right=40] (0.75, 6.5) to (0, 5.25);
            \draw [style=Black Edge, bend left=40] (3.25, 6.5) to (4, 5.25);
            \draw [style=Black Edge, bend right=40] (0.75, 9.5) to (0, 8.25);
            \draw [style=Black Edge, bend left=40] (3.25, 9.5) to (4, 8.25);
            \draw [style=Black Edge, bend right=40] (0.75, 11) to (0, 9.75);
            \draw [style=Black Edge, bend left=40] (3.25, 11) to (4, 9.75);
            \draw [style=Black Edge] (0.75, 2) to (3.25, 2);
            \draw [style=Black Edge] (3.25, 0) to (0.75, 0);
            \draw [style=Black Edge] (0.75, 3.5) to (3.25, 3.5);
            \draw [style=Black Edge] (0, 2.25) to (0, 1);
            \draw [style=Black Edge] (4, 2.25) to (4, 1);
            \draw [style=Black Edge] (0.75, 5) to (3.25, 5);
            \draw [style=Black Edge] (4, 2.25) to (4, 3.75);
            \draw [style=Black Edge] (0, 3.75) to (0, 2.25);
            \draw [style=Black Edge] (0.75, 6.5) to (3.25, 6.5);
            \draw [style=Black Edge] (4, 5.25) to (4, 3.75);
            \draw [style=Black Edge] (0, 5.25) to (0, 3.75);
            \draw [style=Black Edge] (0.75, 9.5) to (3.25, 9.5);
            \draw [style=Black Edge] (4, 8.25) to (4, 5.25);
            \draw [style=Black Edge] (0, 8.25) to (0, 5.25);
            \draw [style=Black Edge] (0.75, 11) to (3.25, 11);
            \draw [style=Black Edge] (0, 9.75) to (0, 8.25);
            \draw [style=Black Edge] (4, 9.75) to (4, 8.25);
            \draw [style=Black Edge] (8) to (17);
            \draw [style=Black Edge] (17) to (10);
            \draw [style=Black Edge] (10) to (15);
            \draw [style=Black Edge] (15) to (16);
            \draw [style=Black Edge] (16) to (13);
            \draw [style=Black Edge] (16) to (14);
            \draw [style=Black Edge] (15) to (17);
            \draw [style=Black Edge] (27) to (22);
            \draw [style=Black Edge] (22) to (24);
            \draw [style=Black Edge] (24) to (23);
            \draw [style=Black Edge] (23) to (26);
            \draw [style=Black Edge] (23) to (25);
            \draw [style=Black Edge] (25) to (22);
            \draw [style=Black Edge] (35) to (32);
            \draw [style=Black Edge] (32) to (33);
            \draw [style=Black Edge] (33) to (34);
            \draw [style=Black Edge] (34) to (32);
            \draw [style=Black Edge] (40) to (44);
            \draw [style=Black Edge] (40) to (41);
            \draw [style=Black Edge] (40) to (43);
            \draw [style=Black Edge] (40) to (42);
            \draw [style=Black Edge] (86) to (84);
            \draw [style=Black Edge] (83) to (86);
            \draw [style=Black Edge] (86) to (82);
            \draw [style=Black Edge] (84) to (85);
            \draw [style=Black Edge] (85) to (87);
            \draw [style=Black Edge] (88) to (87);
            \draw [style=Black Edge] (87) to (89);
            \draw [style=Black Edge] (94) to (92);
            \draw [style=Black Edge] (91) to (94);
            \draw [style=Black Edge] (94) to (90);
            \draw [style=Black Edge] (96) to (95);
            \draw [style=Black Edge] (95) to (97);
            \draw [style=Black Edge] (92) to (95);
            \draw [style=Directed Edge, bend right=45, looseness=1.25] (98) to (99);
            \draw [style=Black Edge, in=-60, out=165] (-0.75, -0.25) to (-1.5, 5.5);
            \draw [style=Black Edge, in=60, out=-165] (100) to (-1.5, 5.5);
        \end{tikzpicture}
    }
\end{figure}
\begin{abstract}
    Phasor Measurement Units (PMUs) are placed at strategic vertices in an electrical power network to monitor the flow of power.
    Determining the minimum number and optimal placement of PMUs is modeled by the graph theoretic process called \emph{Power Domination}.
    This paper describes the \emph{Power Domination Toolbox (PDT)}, which efficiently identifies a minimum number of PMU locations that monitor the entire network.
    The PDT leverages graph theoretic literature to reduce the complexity of determining optimal PMU placements by: reducing the order of the graph (contraction), leveraging zero forcing forts, sorting the remaining solution space, and parallel computing.
    The PDT is a drop-in replacement of the current state-of-the-art exhaustive search algorithm in Python and maintains compatibility with SageMath.
    The PDT can identify minimum PMU placements for graphs with hundreds of vertices on personal computers and can analyze larger graphs on high performance computers.
    The PDT affords users the ability to investigate power domination on graphs previously considered infeasible due to the number of vertices resulting in a prohibitively long run-time.
\end{abstract}
Keywords - optimal sensor placement, phasor measurement units, graph methods, power domination.

\section{Introduction}
    In 2003, a blackout in the Ohio power grid cascaded through a large portion of the Northeastern United States and Canada \cite{Force2004}.
    Blackouts this extensive can be mitigated by monitoring the power grid with phasor measurement units (PMUs) and acting quickly on the information they provide.
    PMUs use conservation of energy laws to observe phasor measures at distant locations in addition to the directly connected power lines.
    One goal for power grid planners is to maximize grid coverage while minimizing the number of installed PMUs to minimize the cost of grid maintenance.
    A minimum set of locations to install PMUs in order to monitor the entire network is called a minimum power dominating set (PDS).

    Currently, some algorithms that are used to find optimal PMU placement locations include: genetic algorithm, particle swarm optimization, tabu search, greedy algorithm, integer linear programming, integer quadratic programming, simulated annealing, hybrid algorithm, exhaustive search, depth-first search, and minimum spanning tree \cite{Baba2021}.
    Recently, Hicks and Smith \cite{Smith2020} have implemented integer linear programming methods in Gurobi to achieve results on large graphs.
    To build the integer linear programming constraints, Hicks and Smith created a restricted problem where PMUs are placed at certain locations and prohibited from being placed at other locations.
    This solution, however, is no longer restricted to the integers and may be slightly greater from the minimum number of PMUs required. 
    To find an exact solution, graph theorists commonly perform power domination via a combination of zero forcing code maintained by Jephian Lin and power domination code developed by Brian Wissman \cite{Benson2018,Hogben2022,Anderson2023,Benson2018a,Bjorkman2022}.
    This method is an exhaustive search for a minimum PDS, and we will call this the \emph{JL-BW algorithm}.
    
    The PDT serves as a drop-in replacement and extension for the JL-BW algorithm written in Python.
    We implement pre-processing techniques used in graph theoretic proofs that have previously not been implemented in software, including: graph contraction to reduce propagation time, leveraging zero forcing forts to restrict the solution space, and assigning a qualitative score to sets in the remaining solution space to determine minimum PDSs more efficiently.
    Parallel compute methodologies are also implemented to fully leverage computational resources.
    To demonstrate the PDT's utility: on random graphs with 120 vertices we see an average 19 times run-time improvement over the JL-BW algorithm.
    The amount of improvement is dependent on the graph structure, and is repeatable on random graphs of other sizes as well as standard test networks.

    This paper will describe the PDT.
    Section~\ref{sec:prelim} will provide the graph theoretic definitions and terminology required in this paper.
    Section~\ref{sec:powdom} will define the power domination algorithm and related concepts.
    Section~\ref{sec:optimization} will describe the process that the PDT uses to find a minimum PDS with brief discussions on the run-time for each step.
    Section~\ref{sec:survey} is a survey of run-time analysis, demonstrating the efficiency of the PDT over the JL-BW algorithm.
    Finally, section~\ref{sec:concluding} gives examples of interfacing with the PDT.

\section{Graph Theory}\label{sec:prelim}
    An electrical power grid can be represented as a \emph{graph} $G$, which consists of two sets: a set of \emph{vertices} (busses), $V(G)$, and a set of unordered pairs of vertices called \emph{edges} (transmission lines), $E(G)$, usually written as $v_1v_2$ for vertices $v_1$ and $v_2$.
    An edge is \emph{incident} to the vertices it contains.
    Two vertices are \emph{adjacent} (neighbors) if there exists an edge between them.
    The degree of a vertex $v$ is the number of vertices adjacent to $v$ and is written $deg_G(v)$.
    When the graph is understood, the subscript is omitted.
    A vertex with degree one is called a \emph{leaf} and a vertex with degree zero is called an \emph{isolated vertex}.
    
    For a graph $G$, a structure $H$ within $G$ is called a \emph{subgraph}, written $H \subseteq G$, when $H$ is a graph with $V(H) \subseteq V(G)$ and $E(H) \subseteq E(G)$.
    A set $A \subseteq V(G)$ generates a \emph{vertex-induced subgraph}, $G[A]$, where $V(G[A]) = A$ and $E(G[A]) = \{ xy \in E(G) : x,y \in A\}$.
    
    A \emph{path} in a graph $G$ has vertices $\left\{v_1, v_2,..., v_n\right\}\subseteq V(G)$, so that $\{v_iv_{i+1}:1\leq i\leq n-1\}\subseteq E(G)$, and is usually written as $v_1v_2...v_n$.
    The \emph{length} of a path is the number of edges it contains.
    A graph $G$ is \emph{connected} if there exists at least one path between any two distinct vertices.
    For a vertex $x$ in $V(G)$, if $G$ is a connected graph but $G[V(G) \setminus \{x\}]$ is not connected, then $x$ is a \emph{cut vertex}.

    For a graph $G$ and edge $xy \in E(G)$, an \emph{edge contraction} (contraction) on $xy$ adds a new vertex $z$ to $G$ such that $z$ is adjacent to any vertex adjacent to either $x$ or $y$ and removes $x$ and $y$ from $G$ as well as any edges incident to $x$ or $y$.

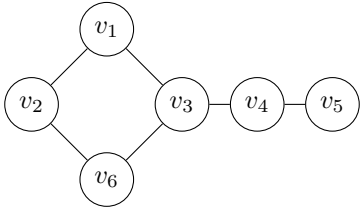
\begin{figure}
        \centering
        \begin{tikzpicture}[every node/.style={draw=black,circle}, minimum size=2em]
            \node (0) at (0, 1) {$v_1$};
            \node (1) at (-1, 0) {$v_2$};
            \node (2) at (0, -1) {$v_6$};
            \node (3) at (1, 0) {$v_3$};
            \node (4) at (2, 0) {$v_4$};
            \node (5) at (3, 0) {$v_5$};
            \draw (1) to (0);
            \draw (0) to (3);
            \draw (3) to (2);
            \draw (2) to (1);
            \draw (3) to (4);
            \draw (4) to (5);
        \end{tikzpicture}
        \caption{The tadpole graph}
        \label{Figure:PendantGraph1}
    \end{figure}

    For a graph $G$ with subset $X \subseteq V(G)$, the \emph{entrance of $X$} in $G$, written $\supp{X}$, is the set of vertices not in $X$ but adjacent to at least one vertex in $X$.
    A \emph{fort}, $F$, in a graph $G$ is a non-empty subset of vertices for which no vertex in $\supp{F}$ is adjacent to exactly one vertex in $F$ \cite{Fast2018}.
    A \emph{terminal fort} with corresponding cut vertex $v$, denoted $F_v$, is a fort in which $\supp{F}=\{v\}$.
    By way of example, the sets $\{v_1, v_2, v_6\}$ and $\{v_1, v_6\}$ are both forts of the tadpole graph in Figure~\ref{Figure:PendantGraph1} with supports $\{v_3\}$ and $\{v_2, v_3\}$ respectively.
    The set $\{v_1, v_2, v_6\}$ is a terminal fort with corresponding cut vertex $v_3$.
    
\section{Power Domination}\label{sec:powdom}
    In the power domination algorithm outlined by Haynes et. al \cite{Haynes2002}, vertices are either \emph{unobserved} or \emph{observed} by a PMU.
    The \emph{power domination algorithm} with input $S \subseteq V(G)$ is as follows:
    \begin{enumerate}
        \item
            (\emph{Domination Step}) Each vertex in $S$, or adjacent to a vertex in $S$, is observed.
        \item
            (\emph{Zero Forcing Step}) While there exists an observed vertex adjacent to exactly one unobserved vertex, the unobserved vertex becomes observed.
    \end{enumerate}
    For a graph $G$ and subset $S \subseteq V(G)$, we identify the set of observed vertices in the graph $G$ after applying the power domination algorithm as $\observed{G}{S}$.
    A \emph{power dominating set (PDS)} is any subset $S\subseteq V(G)$ where $\observed{G}{S}=V(G)$.
    The \emph{power domination number} of a graph $G$ is the cardinality of a minimum PDS, written as $\pdn{G}$.
    Note, there may be many minimum PDSs for a particular graph.

    We also leverage the following terminology from \cite{Bozeman2019}: for a graph $G$ and subset $X \subseteq V(G)$, a \emph{power dominating set subject to $X$} is any subset $S \subseteq V(G)$ containing $X$ that is also a power dominating set.
    In a graph with maximum degree at least 3, a minimum power dominating set can be chosen in which each vertex has degree at least $3$ \cite[Observation~4]{hhhh02} .
    In a restricted power domination problem, certain vertices are already observed by the existing PMUs and may have no unobserved neighbors.
    This means that some vertices provide no additional observations when a PMU is placed on them.
    Combining this with a restriction to vertices with degree at least 3, we define \emph{active vertices} to be vertices that have degree at least 3 and have unobserved neighbors with respect to the restricted power domination problem subject to $X$.
    
\section{Optimizations}\label{sec:optimization}
    The PDT implements optimizations to the power domination process in four ways: 
    \begin{enumerate}
    	\item
    		Contract the input graph.
    	\item
            Leverage zero forcing forts.
    	\item
    		Sort the solution space.

    \end{enumerate}
    Steps 1-3 are done as pre-processing steps and are followed by:

    \begin{enumerate}
    	\item[4.]
    		Distribute the search for a minimum PDS across parallel compute resources.
    \end{enumerate}
    
    The following sections describe each of these optimizations in turn.
    Additionally, we will provide short discussions in each subsection on the run-time impact of implementing the given optimization.
    The discussion centers around a collection of graphs we provide alongside the PDT, and describe in more detail in Section~\ref{sec:carlo}.
    For the purposes of these discussions, this data set is a collection of 600 random graphs with 100 graphs each on 20, 40, 60, 80, 100, and 120 vertices.
    
    \subsection{Contracting the Graph}\label{Section:shrink}
            Due to the nature of the zero forcing step, we have an opportunity to reduce propagation time via contraction.
            Paths on vertices with degree less than 3 are contracted via Algorithm~\ref{alg:contraction}.
            We will then show that this contraction yields a graph with minimum PDSs that are also minimum PDSs of the input graph.
            \begin{algorithm}
                \DontPrintSemicolon
                \KwInput{A graph $G$}
                \KwOutput{A contracted graph $G'$}
                \For{$H$ a connected component of $G\Big[ \{ v \in V(G) : \deg(v) < 3 \} \Big]$}{
                    \If{$H$ is a path terminating in two leaves, or $H$ is a cycle}{
                        Contract the path or cycle in $G$ corresponding to $H$ to an isolated vertex in $G$.
                    }
                    \If{$H$ is a path terminating in vertices adjacent in $G$ to distinct vertices $x, y$ with $\deg_G(x) \geq 3$ and $\deg_G(y) \geq 3$}{
                        Contract the path in $G$ corresponding to $H$ to a single degree 2 vertex.
                    }
                    \If{$H$ is a path terminating in a leaf and a vertex adjacent in $G$ to a vertex $x$ with $\deg_G(x) \geq 3$}{
                        Contract the path in $G$ corresponding to $H$ to a leaf.
                    }
                    \If{$H$ is a path terminating in vertices adjacent in $G$ to a single vertex $x$ with $\deg_G(x) \geq 3$}{
                        Contract the path in $G$ corresponding to $H$ to a pair of adjacent degree 2 vertices.
                    }
                }
                \caption{Graph Contraction Algorithm}
                \label{alg:contraction}
            \end{algorithm}
            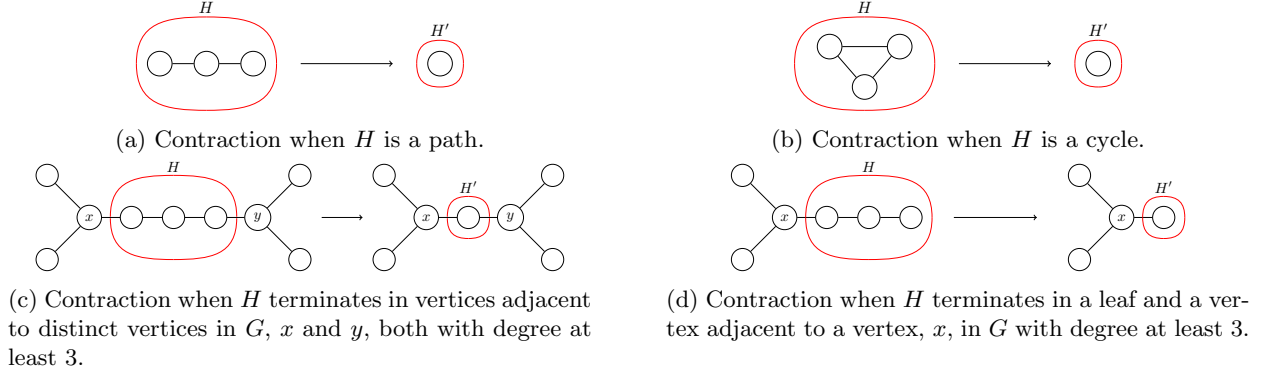
\begin{figure}
                \centering
                \begin{subfigure}{0.47\textwidth}
                    \centering
                    \resizebox{!}{1.5cm}{
                        \begin{tikzpicture}
                            \node [style=Black Outline Node] (0) at (-3, 0) {};
                            \node [style=Black Outline Node] (1) at (-4, 0) {};
                            \node [style=Black Outline Node] (3) at (-2, 0) {};
                            \node [style=Text Node] (13) at (-3, 1.25) {$H$};
                            \node [style=Black Outline Node] (14) at (2, 0) {};
                            \node [style=Text Node] (27) at (2, 0.75) {$H'$};
                            \node [style=Meta Rect] (28) at (-1, 0) {};
                            \node [style=Meta Rect] (29) at (1, 0) {};
                            \draw [style=Black Edge] (1) to (0);
                            \draw [style=Black Edge] (0) to (3);
                            \draw [style=Red Edge, in=90, out=180, looseness=1.25] (-3, 1) to (-4.5, 0);
                            \draw [style=Red Edge, in=-180, out=-90, looseness=1.25] (-4.5, 0) to (-3, -1);
                            \draw [style=Red Edge, in=-90, out=0, looseness=1.25] (-3, -1) to (-1.5, 0);
                            \draw [style=Red Edge, in=360, out=90, looseness=1.25] (-1.5, 0) to (-3, 1);
                            \draw [style=Red Edge, in=90, out=-180, looseness=1.25] (2, 0.5) to (1.5, 0);
                            \draw [style=Red Edge, in=-180, out=-90, looseness=1.25] (1.5, 0) to (2, -0.5);
                            \draw [style=Red Edge, in=-90, out=0, looseness=1.25] (2, -0.5) to (2.5, 0);
                            \draw [style=Red Edge, in=360, out=90, looseness=1.25] (2.5, 0) to (2, 0.5);
                            \draw [style=Directed Edge] (28) to (29);
                        \end{tikzpicture}
                    }
                    \caption{Contraction when $H$ is a path.}
                    \label{subfig:contraction3}
                \end{subfigure}
                \hfill
                \begin{subfigure}{0.47\textwidth}
                    \centering
                    \resizebox{!}{1.5cm}{
                        \begin{tikzpicture}
                            \node [style=Black Outline Node] (0) at (-3, -0.5) {};
                            \node [style=Black Outline Node] (1) at (-3.75, 0.365) {};
                            \node [style=Black Outline Node] (3) at (-2.25, 0.365) {};
                            \node [style=Text Node] (13) at (-3, 1.25) {$H$};
                            \node [style=Black Outline Node] (14) at (2, 0) {};
                            \node [style=Text Node] (27) at (2, 0.75) {$H'$};
                            \node [style=Meta Rect] (28) at (-1, 0) {};
                            \node [style=Meta Rect] (29) at (1, 0) {};
                            \draw [style=Black Edge] (1) to (0);
                            \draw [style=Black Edge] (0) to (3);
                            \draw [style=Black Edge] (1) to (3);
                            \draw [style=Red Edge, in=90, out=180, looseness=1.25] (-3, 1) to (-4.5, 0);
                            \draw [style=Red Edge, in=-180, out=-90, looseness=1.25] (-4.5, 0) to (-3, -1);
                            \draw [style=Red Edge, in=-90, out=0, looseness=1.25] (-3, -1) to (-1.5, 0);
                            \draw [style=Red Edge, in=360, out=90, looseness=1.25] (-1.5, 0) to (-3, 1);
                            \draw [style=Red Edge, in=90, out=-180, looseness=1.25] (2, 0.5) to (1.5, 0);
                            \draw [style=Red Edge, in=-180, out=-90, looseness=1.25] (1.5, 0) to (2, -0.5);
                            \draw [style=Red Edge, in=-90, out=0, looseness=1.25] (2, -0.5) to (2.5, 0);
                            \draw [style=Red Edge, in=360, out=90, looseness=1.25] (2.5, 0) to (2, 0.5);
                            \draw [style=Directed Edge] (28) to (29);
                        \end{tikzpicture}
                    }
                    \caption{Contraction when $H$ is a cycle.}
                    \label{subfig:contraction4}
                \end{subfigure}
                \newline
                \begin{subfigure}{0.47\textwidth}
                    \centering
                    \resizebox{!}{1.5cm}{
                    \begin{tikzpicture}
                        \node [style=Black Outline Node] (0) at (-4, 0) {};
                        \node [style=Black Outline Node] (1) at (-5, 0) {};
                        \node [style=Black Outline Node] (2) at (-6, 0) {$x$};
                        \node [style=Black Outline Node] (3) at (-3, 0) {};
                        \node [style=Black Outline Node] (4) at (-2, 0) {$y$};
                        \node [style=Black Outline Node] (5) at (-1, 1) {};
                        \node [style=Black Outline Node] (6) at (-1, -1) {};
                        \node [style=Black Outline Node] (7) at (-7, 1) {};
                        \node [style=Black Outline Node] (8) at (-7, -1) {};
                        \node [style=Meta Rect] (10) at (-5.5, 0) {};
                        \node [style=Meta Rect] (12) at (-2.5, 0) {};
                        \node [style=Text Node] (13) at (-4, 1.25) {$H$};
                        \node [style=Black Outline Node] (14) at (3, 0) {};
                        \node [style=Black Outline Node] (16) at (2, 0) {$x$};
                        \node [style=Black Outline Node] (18) at (4, 0) {$y$};
                        \node [style=Black Outline Node] (19) at (5, 1) {};
                        \node [style=Black Outline Node] (20) at (5, -1) {};
                        \node [style=Black Outline Node] (21) at (1, 1) {};
                        \node [style=Black Outline Node] (22) at (1, -1) {};
                        \node [style=Meta Rect] (24) at (2.5, 0) {};
                        \node [style=Meta Rect] (26) at (3.5, 0) {};
                        \node [style=Text Node] (27) at (3, 0.75) {$H'$};
                        \node [style=Meta Rect] (28) at (-0.5, 0) {};
                        \node [style=Meta Rect] (29) at (0.5, 0) {};
                        \draw [style=Black Edge] (7) to (2);
                        \draw [style=Black Edge] (2) to (8);
                        \draw [style=Black Edge] (2) to (1);
                        \draw [style=Black Edge] (1) to (0);
                        \draw [style=Black Edge] (0) to (3);
                        \draw [style=Black Edge] (3) to (4);
                        \draw [style=Black Edge] (4) to (5);
                        \draw [style=Black Edge] (4) to (6);
                        \draw [style=Red Edge, in=90, out=180, looseness=1.25] (-4, 1) to (10);
                        \draw [style=Red Edge, in=-180, out=-90, looseness=1.25] (10) to (-4, -1);
                        \draw [style=Red Edge, in=-90, out=0, looseness=1.25] (-4, -1) to (12);
                        \draw [style=Red Edge, in=360, out=90, looseness=1.25] (12) to (-4, 1);
                        \draw [style=Black Edge] (21) to (16);
                        \draw [style=Black Edge] (16) to (22);
                        \draw [style=Black Edge] (18) to (19);
                        \draw [style=Black Edge] (18) to (20);
                        \draw [style=Red Edge, in=90, out=-180, looseness=1.25] (3, 0.5) to (24);
                        \draw [style=Red Edge, in=-180, out=-90, looseness=1.25] (24) to (3, -0.5);
                        \draw [style=Red Edge, in=-90, out=0, looseness=1.25] (3, -0.5) to (26);
                        \draw [style=Red Edge, in=360, out=90, looseness=1.25] (26) to (3, 0.5);
                        \draw [style=Directed Edge] (28) to (29);
                        \draw [style=Black Edge] (16) to (14);
                        \draw [style=Black Edge] (14) to (18);
                    \end{tikzpicture}
                    }
                    \caption{Contraction when $H$ terminates in vertices adjacent to distinct vertices in $G$, $x$ and $y$, both with degree at least 3.}
                    \label{subfig:contraction1}
                \end{subfigure}
                \hfill
                \begin{subfigure}{0.47\textwidth}
                    \centering
                    \resizebox{!}{1.5cm}{
                        \begin{tikzpicture}
                            \node [style=Black Outline Node] (0) at (-3, 0) {};
                            \node [style=Black Outline Node] (1) at (-4, 0) {};
                            \node [style=Black Outline Node] (2) at (-5, 0) {$x$};
                            \node [style=Black Outline Node] (3) at (-2, 0) {};
                            \node [style=Black Outline Node] (7) at (-6, 1) {};
                            \node [style=Black Outline Node] (8) at (-6, -1) {};
                            \node [style=Text Node] (13) at (-3, 1.25) {$H$};
                            \node [style=Black Outline Node] (14) at (4, 0) {};
                            \node [style=Black Outline Node] (18) at (3, 0) {$x$};
                            \node [style=Black Outline Node] (19) at (2, 1) {};
                            \node [style=Black Outline Node] (20) at (2, -1) {};
                            \node [style=Text Node] (27) at (4, 0.75) {$H'$};
                            \node [style=Meta Rect] (28) at (-1, 0) {};
                            \node [style=Meta Rect] (29) at (1, 0) {};
                            \draw [style=Black Edge] (7) to (2);
                            \draw [style=Black Edge] (2) to (8);
                            \draw [style=Black Edge] (2) to (1);
                            \draw [style=Black Edge] (1) to (0);
                            \draw [style=Black Edge] (0) to (3);
                            \draw [style=Black Edge] (18) to (19);
                            \draw [style=Red Edge, in=90, out=180, looseness=1.25] (-3, 1) to (-4.5, 0);
                            \draw [style=Red Edge, in=-180, out=-90, looseness=1.25] (-4.5, 0) to (-3, -1);
                            \draw [style=Red Edge, in=-90, out=0, looseness=1.25] (-3, -1) to (-1.5, 0);
                            \draw [style=Red Edge, in=360, out=90, looseness=1.25] (-1.5, 0) to (-3, 1);
                            \draw [style=Black Edge] (18) to (20);
                            \draw [style=Black Edge] (14) to (18);
                            \draw [style=Red Edge, in=90, out=-180, looseness=1.25] (4, 0.5) to (3.5, 0);
                            \draw [style=Red Edge, in=-180, out=-90, looseness=1.25] (3.5, 0) to (4, -0.5);
                            \draw [style=Red Edge, in=-90, out=0, looseness=1.25] (4, -0.5) to (4.5, 0);
                            \draw [style=Red Edge, in=360, out=90, looseness=1.25] (4.5, 0) to (4, 0.5);
                            \draw [style=Directed Edge] (28) to (29);
                        \end{tikzpicture}
                    }
                    \caption{Contraction when $H$ terminates in a leaf and a vertex adjacent to a vertex, $x$, in $G$ with degree at least 3. \phantom{this is phantom text}}
                    \label{subfig:contraction2}
                \end{subfigure}
                \caption{Cases where Algorithm \ref{alg:contraction} contracts $H$ to a single vertex.}
                \label{fig:singletoncontractions}
            \end{figure}
            \begin{figure}
                \centering
                \resizebox{!}{2cm}{
                    \begin{tikzpicture}
                        \node [style=Black Outline Node] (0) at (-4.5, 0) {$x$};
                        \node [style=Black Outline Node] (1) at (-3.5, 1) {};
                        \node [style=Black Outline Node] (4) at (-3.5, -1) {};
                        \node [style=Black Outline Node] (6) at (-5.5, 0) {};
                        \node [style=Black Outline Node] (7) at (-2.5, 0) {};
                        \node [style=Text Node] (13) at (-1.75, 0) {$H$};
                        \node [style=Black Outline Node] (16) at (2, 0) {};
                        \node [style=Black Outline Node] (17) at (3, 0) {$x$};
                        \node [style=Black Outline Node] (18) at (4, 1) {};
                        \node [style=Black Outline Node] (19) at (4, -1) {};
                        \node [style=Text Node] (24) at (4.75, 0) {$H'$};
                        \node [style=Meta Rect] (14) at (-1, 0) {};
                        \node [style=Meta Rect] (15) at (1, 0) {};
                        \draw [style=Black Edge] (6) to (0);
                        \draw [style=Black Edge] (0) to (1);
                        \draw [style=Black Edge] (1) to (7);
                        \draw [style=Black Edge] (7) to (4);
                        \draw [style=Black Edge] (4) to (0);
                        \draw [style=Red Edge, bend left=90] (-4, 1) to (-2, 1);
                        \draw [style=Red Edge] (-4, 1) to (-4, -1);
                        \draw [style=Red Edge, bend right=90] (-4, -1) to (-2, -1);
                        \draw [style=Red Edge] (-2, 1) to (-2, -1);
                        \draw [style=Black Edge] (17) to (18);
                        \draw [style=Black Edge] (18) to (19);
                        \draw [style=Black Edge] (19) to (17);
                        \draw [style=Black Edge] (17) to (16);
                        \draw [style=Red Edge, bend left=90, looseness=1.5] (3.5, 1) to (4.5, 1);
                        \draw [style=Red Edge] (4.5, 1) to (4.5, -1);
                        \draw [style=Red Edge, bend left=90, looseness=1.5] (4.5, -1) to (3.5, -1);
                        \draw [style=Red Edge] (3.5, -1) to (3.5, 1);
                        \draw [style=Directed Edge] (14) to (15);
                    \end{tikzpicture}
                }
                \caption{Case where Algorithm \ref{alg:contraction} contracts $H$ to a pair of adjacent vertices: when $H$ terminates in vertices adjacent to a unique vertex in $G$ with degree at least 3.}
                \label{fig:adjacentcontractions}
            \end{figure}
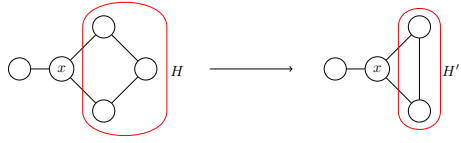

            Algorithm \ref{alg:contraction} contracts subgraphs of the input graph to either a single vertex, or a pair of adjacent vertices.
            Figure \ref{fig:singletoncontractions} provides examples of when Algorithm \ref{alg:contraction} contracts $H$ to a single vertex which corresponds to the conditional statements on lines 2, 4, and 6.
            Figure \ref{fig:adjacentcontractions} provides the example of when Algorithm \ref{alg:contraction} contracts $H$ to a pair of adjacent vertices which corresponds to the conditional statement on line 8.

            A minimum power dominating set for the contracted graph resulting from Algorithm~\ref{alg:contraction} corresponds to a minimum power dominating set of the input graph.

            \begin{theorem}
                Let $G$ be a graph and let $G'$ be the result of contracting $G$ via Algorithm~\ref{alg:contraction}.
                A minimum power dominating set of $G'$ corresponds to a minimum power dominating set of $G$. Therefore, $\pdn{G}=\pdn{G'}$.
            \end{theorem}

            \begin{proof}
                We will begin with a minimum power dominating set $S$ of $G$ and demonstrate a corresponding minimum power dominating set of $G'$.
                For any vertex $v\in S$, if $v\not\in V(G')$, then $v$ is a degree 1 or 2 vertex in $G$ that was contracted to create $G'$.
                Note that as $S$ is minimum, there cannot be 2 adjacent such vertices in $S$.
                Replace each vertex $v$ with its corresponding contracted vertex and call this set $S'$.
                The power domination process on $G'$ with initial set $S'$ proceeds analogously to the power domination process on $G$ with initial set $S$, with removed observations along the contracted paths.
                
                Next, consider a power dominating set $S'$ of $G'$.
                For any vertex $x\in S'$, if $x\not\in V(G)$, then $x$ is the result of contracting degree 1 or 2 vertices in $G$.
                Replace each such $x$ with one of the corresponding vertices and call the resulting set $S$.
                The power domination process on $G$ with initial set $S$ proceeds analogously to the power domination process on $G'$ with initial set $S'$, with added observations along the non-contracted paths.

                We have shown corresponding minimum power dominating sets of the same size for $G$ and $G'$ and so $\pdn{G}=\pdn{G'}$.
            \end{proof}

            Contracting the graph is done in linear time, while the propagation steps it eliminates from the power domination process is potentially exponential.
            This is a novel method of the PDT that dramatically improves run-time on graphs that demonstrate substructures of long chains compared to the JL-BW algorithm.
            Table~\ref{tab:contractiontimes} outlines the average run-time to contract random graphs of varying size.
            \begin{table}[ht]
                \centering
                \resizebox{\textwidth}{!}{
                    \begin{tabular}{|c|c|c|c|c|c|c|}
                        \hline
                        $|V(G)|$                 & 20                   & 40                   & 60                   & 80                    & 100                  & 120                  \\ \hline
                        Time to   determine $G'$ & $2.320\times10^{-4}$ & $4.357\times10^{-4}$ & $6.440\times10^{-4}$ & $6.907\times10^{-4}$ & $1.251\times10^{-3}$ & $1.661\times10^{-3}$ \\ \hline
                    \end{tabular}
                }
                \caption{Average time (in seconds) the PDT uses to calculate $G'$ for graphs of varying sizes.}
                \label{tab:contractiontimes}
            \end{table}
        
    \subsection{Leveraging Forts}\label{Section:guaranteed}
        Forts were first utilized in zero forcing by Fast and Hicks in \cite{Fast2018}. 
        They were generalized to power domination Bozeman et al. in \cite{Bozeman2019}.
        Hicks and Smith exploited forts to find minimum power dominating sets in their integer linear program method \cite{Smith2020}.
        We will also use forts in order to find a minimum power dominating set while considering significantly fewer cases than brute force methods.

        \begin{proposition}[\hspace{-1pt}{\cite[Proposition~4.3]{Bozeman2019}}]\label{prop:bozemanforts}
            Let $G$ be a graph and $F$ be any fort of $G$. If $S$ is a power dominating set of $G$, then $S\cap N[F]\neq \emptyset.$
        \end{proposition}

        Thus, any power dominating set must intersect with the closed neighborhood of every fort.
        In special cases, we can be more particular about which vertices from the closed neighborhood of a fort are in some minimum power dominating set. 
        Note that $N[F] = F\cup \supp{F}$.
    
        \begin{theorem}\label{thm:unifier}         
            There exists a minimum power dominating set $S$ such that $S\cap F = \emptyset$ and $S\cap \supp{F} \neq \emptyset$ for every fort $F$ satisfying $F\cup \supp{F} \subseteq \observed{G}{\{v\}}$ for all $v\in \supp{F}$.
        \end{theorem}
        \begin{proof}
            Let $S$ be a minimum power dominating set of $G$. Let $F$ be a fort satisfying $F\cup \supp{F} \subseteq \observed{G}{\{v\}}$ for all $v\in \supp{F}$.

            If $S\cap F\neq \emptyset$, then there is some vertex $w\in S \cap F$.
            We construct $S'= (S\setminus \{w\})\cup\{v\}$ for any $v\in \supp{F}$.
            As $F\cup \supp{F} \subseteq\observed{G}{\{v\}}$, and any neighbor of $v$ is either in $F$ or in $\supp{F}$, any vertex observed as  a result of $v$ is also observed as a result of $w$. Thus $S'$ is a power dominating set of the same size as $S$.
    
            If $S\cap \supp{F} = \emptyset$ then by Proposition \ref{prop:bozemanforts}, there is some vertex $w \in S \cap F$, and the first case is recovered.
        \end{proof}

        Theorem \ref{thm:unifier} can be used to determine a set of vertices that are in \emph{some} minimum power dominating set of a given graph.
        Determining all such vertices may be time-intensive, and so we focus on 2 cases: terminal forts and forts associated with induced $C_4$ subgraphs in the contracted graph. 

        \begin{observation}
        For any graph $G$ with corresponding contracted graph $G'$ from Algorithm~\ref{alg:contraction}, there exists a minimum power dominating set $S$ of $G'$ such that
        \begin{enumerate}
            \item
                $S\cap F = \emptyset$ and $v\in S$ for every terminal fort $F_v$ satisfying $F_v\cup \{v\} \subseteq \observed{G'}{\{v\}}$, and
            \item
                For any zero forcing fort $F = \{y,z\}$ with $\deg_G(y)=\deg_G(z)=2$ so that $F\cup \supp{F}$ is an induced $C_4$ subgraph in $G'$, we have $S\cap F = \emptyset$ and $S\cap \supp{F} \neq \emptyset$.
                Note that in this case, $F$ automatically satisfies $F\cup \supp{F} \subseteq \observed{G'}{\{v\}}$ for all $v\in \supp{F}$.
        \end{enumerate}
        \end{observation}

        That is, for any contracted graph $G'$ from Algorithm~\ref{alg:contraction}, there exists a minimum power dominating set $S$ containing: all isolated vertices, the entrance of every terminal fort, and containing at least one of the 2 entrance vertices of any induced $C_4$ in $G'$ containing 2 non-adjacent vertices with degree 2.
        We will call the entrance vertices of terminal forts \emph{preferred vertices}, and let $\pref{G} = \{ v : F_v \text{ is a terminal fort}\}$.
        Let all \emph{paired entrance vertices} to zero forcing forts corresponding to an induced $C_4$ containing 2 non-adjacent vertices with degree 2 be given as $\jeff{G}=$.

        Preferred vertices and paired entrance vertices give us lower bounds on the power domination number.  
        Let $\Gamma = \Big\{\{\xi\}: \xi \in \pref{G'}\Big\} \cup \Big\{ \xi \in \jeff{G'}\Big\}$ and let $\Phi$ be a minimum set of vertices of $G'$ that intersects non-trivially with each set in $\Gamma$.
        Calculating $\Phi$ explicitly can be done by creating a graph with vertices $V=\Gamma$ and edges $E=\{xy : x\cap y \neq \emptyset\}$ and determining a dominating set. 
        Instead, we approximate $|\Phi|$ by the number of components in this graph, written as $\varphi$.

        \begin{observation}
            For any graph $G$ with $\pdn{G} >1$ and corresponding contraction $G'$ from Algorithm~\ref{alg:contraction},
            \begin{align*}
                \pdn{G} \geq |\Phi| \geq \varphi \geq |\pref{G'}|.
            \end{align*}
        \end{observation}

        The condition $\pdn{G} > 1$ consists of an edge case in which terminal forts overlap, yielding multiple preferred vertices for a graph with power domination number one.
        For example, consider $G$ to be $C_3$ with a leaf added to distinct vertices $a$ and $b$. $\pdn{G} = 1$ but $\pref{G} = \{a, b\}$. This can be shown to occur precisely when $\pdn{G} = 1$.

        In practice, the PDT determines preferred vertices first, as shown in Algorithm~\ref{alg:prefalgorithm}.
        The conditional on line 6 corresponds to catching the edge case where $\pdn{G} = 1$, which is nested in the loop on line 5.
        The loop on line 5 locates terminal forts, including type I and type II forts as defined by Hicks and Smith \cite{Smith2020} and more general terminal forts as shown in Figure~\ref{fig:terminalforts}.
        Paired entrance vertices are located with Algorithm~\ref{alg:c4algorithm}, which leverages fast methods for determining chordless cycles within a graph.
        The run-time to calculate $\pref{G}$ and $\jeff{G}$ is given in Table~\ref{tab:prefjeffruntime}.
        \begin{algorithm}
            \DontPrintSemicolon
            \KwInput{A graph $G$}
            \KwOutput{$\pref{G}$}
            ${\pref{G}} \gets \emptyset$\\            
            \For{$v \in V(G)$ where $v$ is a cut vertex}{
                \If{$v$ is adjacent to at least 2 leaves}{
                    ${\pref{G}} \gets {\text{Pref}(G) \cup {\{v\}}}$
                }
                \For{$H$ a connected component of $G[V(G)\setminus \{v\}]$, and $|V(H)| > 1$}{
                    \If{$H \subseteq \observed{G}{\{v\}}$}{
                        \If{$\pref{G} = \emptyset$}
                        {
                            \If{$V(G) \subseteq \observed{G}{\{v\}}$}
                            {
                                \Return $\{v\}$
                            }
                        }
                        ${\pref{G}} \gets {\text{Pref}(G) \cup {\{v\}}}$
                    }
                }
            }
            \Return $\pref{G}$
            \caption{Algorithm to determine Pref$(G)$}
            \label{alg:prefalgorithm}
        \end{algorithm}
        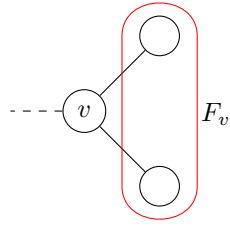
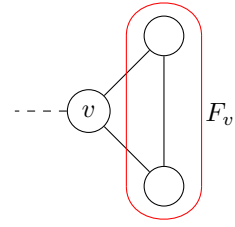
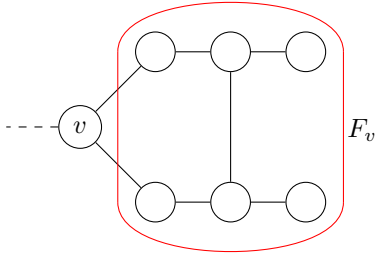
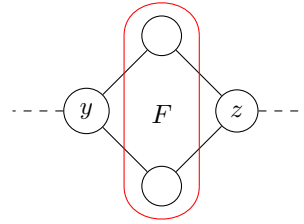
\begin{figure}
            \centering
            \begin{subfigure}{0.49\textwidth}
                \centering
                \begin{tikzpicture}
                    \node [style=Meta Rect] (16) at (2, 0) {};
                    \node [style=Black Outline Node] (17) at (3, 0) {$v$};
                    \node [style=Black Outline Node] (18) at (4, 1) {};
                    \node [style=Black Outline Node] (19) at (4, -1) {};
                    \node [style=Text Node] (24) at (4.75, 0) {$F_v$};
                    \draw [style=Black Edge] (17) to (18);
                    \draw [style=Black Edge] (19) to (17);
                    \draw [style=dashed] (17) to (16);
                    \draw [style=Red Edge, bend left=90, looseness=1.5] (3.5, 1) to (4.5, 1);
                    \draw [style=Red Edge] (4.5, 1) to (4.5, -1);
                    \draw [style=Red Edge, bend left=90, looseness=1.5] (4.5, -1) to (3.5, -1);
                    \draw [style=Red Edge] (3.5, -1) to (3.5, 1);
                \end{tikzpicture}
                \caption{A pair of leaves, $F_v$, indicated by Hicks and Smith as Type I forts \cite{Smith2020}.}
            \end{subfigure}
            \hfill
            \begin{subfigure}{0.49\textwidth}
                \centering
                \begin{tikzpicture}
                    \node [style=Meta Rect] (16) at (2, 0) {};
                    \node [style=Black Outline Node] (17) at (3, 0) {$v$};
                    \node [style=Black Outline Node] (18) at (4, 1) {};
                    \node [style=Black Outline Node] (19) at (4, -1) {};
                    \node [style=Text Node] (24) at (4.75, 0) {$F_v$};
                    \draw [style=Black Edge] (17) to (18);
                    \draw [style=Black Edge] (18) to (19);
                    \draw [style=Black Edge] (19) to (17);
                    \draw [style=dashed] (17) to (16);
                    \draw [style=Red Edge, bend left=90, looseness=1.5] (3.5, 1) to (4.5, 1);
                    \draw [style=Red Edge] (4.5, 1) to (4.5, -1);
                    \draw [style=Red Edge, bend left=90, looseness=1.5] (4.5, -1) to (3.5, -1);
                    \draw [style=Red Edge] (3.5, -1) to (3.5, 1);
                \end{tikzpicture}
                \caption{A terminal $C_3$, $F_v$, indicated by Hicks and Smith as Type II forts \cite{Smith2020}.}
            \end{subfigure}
            \newline
            \begin{subfigure}{0.49\textwidth}
                \centering
                \begin{tikzpicture}
                    \node [style=Meta Rect] (16) at (2, 0) {};
                    \node [style=Black Outline Node] (17) at (3, 0) {$v$};
                    \node [style=Black Outline Node] (18) at (4, 1) {};
                    \node [style=Black Outline Node] (19) at (4, -1) {};
                    \node [style=Black Outline Node] (20) at (5, 1) {};
                    \node [style=Black Outline Node] (21) at (5, -1) {};
                    \node [style=Black Outline Node] (22) at (6, 1) {};
                    \node [style=Black Outline Node] (23) at (6, -1) {};
                    \node [style=Text Node] (24) at (6.75, 0) {$F_v$};
                    \draw [style=Black Edge] (17) to (18);
                    \draw [style=Black Edge] (19) to (17);
                    \draw [style=dashed edge] (17) to (16);
                    \draw [style=Black Edge] (18) to (20);
                    \draw [style=Black Edge] (19) to (21);
                    \draw [style=Black Edge] (20) to (21);
                    \draw [style=Black Edge] (20) to (22);
                    \draw [style=Black Edge] (21) to (23);
                    \draw [style=Red Edge, bend left=90, looseness=0.75] (3.5, 1) to (6.5, 1);
                    \draw [style=Red Edge] (6.5, 1) to (6.5, -1);
                    \draw [style=Red Edge, bend left=90, looseness=0.75] (6.5, -1) to (3.5, -1);
                    \draw [style=Red Edge] (3.5, -1) to (3.5, 1);
                \end{tikzpicture}
                \caption{An example of a terminal fort, $F_v$, not described by Hicks and Smith, that is located by the PDT.}
            \end{subfigure}
            \hfill
            \begin{subfigure}{0.49\textwidth}
                \centering
                \begin{tikzpicture}
                    \node [style=Meta Rect] (16) at (2, 0) {};
                    \node [style=Black Outline Node] (17) at (3, 0) {\resizebox{5pt}{!}{$y$}};
                    \node [style=Black Outline Node] (18) at (4, 1) {};
                    \node [style=Black Outline Node] (19) at (4, -1) {};
                    \node [style=Black Outline Node] (20) at (5, 0) {\resizebox{5pt}{!}{$z$}};
                    \node [style=Meta Rect] (21) at (6, 0) {};
                    \node [style=Text Node] (24) at (4, 0) {$F$};
                    \draw [style=Black Edge] (17) to (18);
                    \draw [style=Black Edge] (19) to (17);
                    \draw [style=dashed] (17) to (16);
                    \draw [style=Black Edge] (18) to (20);
                    \draw [style=Black Edge] (19) to (20);
                    \draw [style=dashed] (20) to (21);            
                    \draw [style=Red Edge, bend left=90, looseness=1.5] (3.5, 1) to (4.5, 1);
                    \draw [style=Red Edge] (4.5, 1) to (4.5, -1);
                    \draw [style=Red Edge, bend left=90, looseness=1.5] (4.5, -1) to (3.5, -1);
                    \draw [style=Red Edge] (3.5, -1) to (3.5, 1);
                \end{tikzpicture}
                \caption{A fort associated with an induced $C_4$, $F$, indicated by Hicks and Smith as Type III forts \cite{Smith2020}.}
            \end{subfigure}
            \caption{Forts located by the PDT indicated in red. Dashed edges represent connection(s) from the entrance vertices to the remainder of the graph.}
            \label{fig:terminalforts}
        \end{figure}
        \begin{algorithm}
            \DontPrintSemicolon
            \KwInput{A graph $G$}
            \KwOutput{$\jeff{G}$}
            $G' \gets {\text{ Algorithm~\ref{alg:contraction}}(G)}$\\
            $A \gets {\{v \in V(G') : \deg_G{v} = 2\}}$\\
            $B \gets N_{G'}[A]$\\
            $\jeff{G} \gets \emptyset$\\            
            \For{$H$ a connected component of $G'[B]$}{
                \If{$C_4 \subseteq H$}{
                    $\jeff{G} \gets {\jeff{G} + {\{v \in V(H) : \deg_G(v) > 2\}}}$
                }
            }
            \caption{Algorithm to determine the entrance of forts associated with induced $C_4$ subgraphs of the contracted graph}
            \label{alg:c4algorithm}
        \end{algorithm}
        \begin{table}
            \resizebox{\textwidth}{!}{%
                \begin{tabular}{|c|c|c|c|c|c|c|}
                \hline
                $|V(G)|$                           & 20                   & 40                   & 60                  & 80                   & 100                  & 120                  \\ \hline
                Time to determine $\pref{G}$ & $4.898\times10^{-4}$ s & $1.234\times10^{-3}$ s & $2.046\times10^{-3}$ s & $2.501\times10^{-3}$ s & $2.703\times10^{-3}$ s & $2.831\times10^{-3}$ s \\ \hline
                Time to determine $\jeff{G}$ & $4.225\times10^{-4}$ s & $9.010\times10^{-4}$ s & $1.385\times10^{-3}$ s & $1.487\times10^{-3}$ s & $1.656\times10^{-3}$ s & $1.905\times10^{-3}$ s \\ \hline
                \end{tabular}%
            }
            \caption{Average run-time (in seconds) the PDT uses to calculate $\pref{G}$ and $\jeff{G}$.}
            \label{tab:prefjeffruntime}
        \end{table}
                
    \subsection{Determining Qualitative Scores}\label{Section:likely}
        The PDT sorts the potential power dominating sets in the solution space to more optimally locate power dominating sets.
        This is done by maximizing the number of vertices observed after considering the restricted power domination problem on $G'$ subject to $\pref{G'}$.
        Define the \emph{qualitative score}, of a vertex $v$ as
        \begin{equation}\label{Equation:likelihood}
            Q(v):=\observed{G}{\pref{G}\cup\{v\}},
        \end{equation}
        Formulating $Q(v)$ in this way affords another opportunity for optimization.
        If $\pdn{G} \neq |\pref{G}|$ and $\max(Q(v):v \in V(G))=|V(G)|$, then $\pdn{G} = |\pref{G}|+1$ and $\pref{G}\cup\{v\}$ is a minimum PDS for any $v \in V(G)$ such that $Q(v) = |V(G)|$.
        
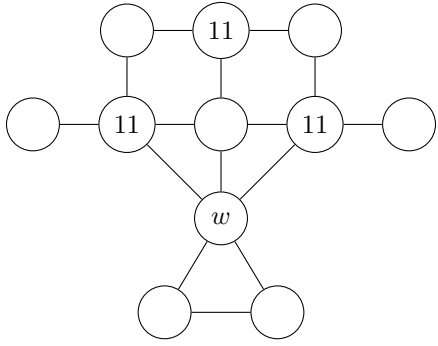
\begin{figure}
            \centering
            \begin{tikzpicture}[every node/.style={draw=black,circle}, minimum size=2em]
                \node (3) at (3, -0.25) {};
                \node (4) at (4.25, -0.25) {11};
                \node (5) at (4.25, 1) {};
                \node (6) at (5.5, 1) {11};
                \node (7) at (5.5, -0.25) {};
                \node (8) at (5.5, -1.5) {$w$};
                \node (9) at (6.75, -0.25) {11};
                \node (10) at (6.75, 1) {};
                \node (11) at (8, -0.25) {};
                \node (13) at (4.75, -2.75) {};
                \node (18) at (6.25, -2.75) {};
                \draw (3) to (4);
                \draw (4) to (8);
                \draw (8) to (13);
                \draw (8) to (18);
                \draw (9) to (8);
                \draw (7) to (8);
                \draw (4) to (7);
                \draw (7) to (9);
                \draw (9) to (11);
                \draw (10) to (9);
                \draw (10) to (6);
                \draw (6) to (7);
                \draw (5) to (4);
                \draw (5) to (6);
                \draw (13) to (18);
            \end{tikzpicture}
            \caption{Zim graph with qualitative score for each non-preferred vertex with degree at least 3 displayed. Vertex $w$ is the only preferred vertex.}
            \label{Figure:ZimGraphRatings}
        \end{figure}

        The PDT iteratively adds vertices to a potential PDS from largest to smallest qualitative score.
        Figure~\ref{Figure:ZimGraphRatings} provides $Q(v)$ for each active vertex in the Zim graph.
        The PDT checks subsets with higher total qualitative score first.
        The JL-BW algorithm does no such pre-processing and checks subsets lexicographically by vertex label.
    
    \subsection{Distribution Across Multiple Threads}\label{Section:parallel}
        Using an exhaustive search algorithm to determine the power domination number of a graph, necessarily, requires determining $\observed{G}{S}$ for sets $S \subseteq V(G)$ with $|S| < \pdn{G}$.
        While the PDT dramatically lowers this number of sets $S$, the number of sets to analyze with the PDT still grows factorially as the power domination number grows.
        Checking each set $S$ can be viewed as an independent process and so the PDT utilizes parallelization.
        
        The multiprocessing library in Python is implemented to distribute the search for a minimum PDS across available CPU resources.
        This parallelization requires computational overhead, so the PDT checks the first 50,000 subsets that may be PDSs on a single compute thread.
        Once the number of subsets to check grows, however, it becomes advantageous to parallelize the search for a minimum PDS with a given number of PMUs.
        The PDT can dynamically determine the number of processes to leverage, or it can be set explicitly to use a predetermined number of processes.
        Common personal computers can facilitate the search for a minimum PDS across approximately 10 processes while allowing the user to continue using the computer for light tasks.

    \subsection{The PDT Algorithm}        
        We now discuss the algorithm that the PDT uses to find a minimum PDS of an input graph $G$.
        The PDT calculates $G'$ and the parameters $\pref{G'}$, $\jeff{G'}$, $\observed{G'}{\pref{G'}}$, and active vertices.
        The PDT then iterates over the connected components of $G'$, restricts the predetermined parameters to the connected component, and calls Algorithm~\ref{alg:connectedalgorithm} with these parameters as additional input to determine a minimum PDS of the component.
        The PDT allows the user to directly call this algorithm as \verb|PDT_minpds_connected|, and determines $G'$, $\pref{G'}$, $\jeff{G'}$, $\observed{G'}{\pref{G'}}$, and active vertices if not provided by the user.
        Let us inspect Algorithm~\ref{alg:connectedalgorithm} as if the user called it directly.
        
        Line 1 is the application of the graph contraction algorithm.
        If there are no vertices with degree more than 2, then any vertex is a PDS, and hence an arbitrary vertex is returned on line 4.

        The PDT then calculates $\pref{G'}$, and $\observed{G}{\pref{G'}}$ and stores these observed vertices in the set $B$ on line 5.
        If $B = V(G)$, then $\pref{G'}$ is a PDS and is then returned on line 8.
        Lines 9 and 10 then locate active vertices for the restricted power domination problem on $G'$ subject to $\pref{G'}$ and determine $Q(v)$ for each active vertex.
        If $\max(Q(v)) = |V(G')|$, then the PDT has located a minimum power dominating set with a PMU on each vertex in $\pref{G'}$ and a vertex with $Q(v)=|V(G')|$.
        The PDT returns this minimum PDS on line 13.

        If a PDS has not been located at this point, $\pdn{G'} \geq |\pref{G'}|+2$.
        The PDT then determines $\jeff{G'}$ and $\varphi$.      
        If $\varphi-|\pref{G}|>2$, more than 2 additional vertices are required to form a PDS.
        These two cases are covered by $i = \max\{2, \varphi-|\pref{G}|\}$ on line 14.
        
        The PDT then begins checking sets of size $|\pref{G'}| + i$ for a minimum PDS.
        Line 17 sets up the combination of additional vertices, $C$, to be added to $\pref{G'}$ to create a subset $S$, which is formed on line 18.
        The PDT enforces that $S$ intersects non-trivially with each pair of entrance vertices on line 19.
        Determining the observed vertices in $\observed{G}{\pref{G'}\cup C}$ on line 20 is given by repeated application of the zero forcing step to the set $\observed{G}{\pref{G}} \cup N[C]$.
        If the resulting set is equal to $V(G')$, then the minimum PDS $S$ is returned on line 22.
        If all $B$ with $|B|=i$ are exhausted, then $i$ is incremented on line 23 and the PDT returns to the loop starting on line 16.
        \begin{algorithm}
            \DontPrintSemicolon
            \KwInput{A \emph{connected} graph $G$}
            \KwOutput{A minimum PDS of $G$}
            ${G' \gets}$ Algorithm~\ref{alg:contraction}$(G)$\\
            \If{$\not \exists v \in V(G)$ with $\deg_G(v) > 2$}
            {
                $S \gets \{v\}$ for any $v \in V(H)$\\
                \Return $S$
            }
            $B \gets \observed{G'}{\pref{G'}}$\\
            \If{$B = V(G')$}
            {
                $S \gets \pref{G'}$\\
                \Return $S$
            }
            $U \gets V(G') \setminus B$\\
            $A \gets {\{v \in V(G') : \deg_{G'}(v) > 2 \text{ and } N[v] \cap U \neq \emptyset\}}$\\
            \If{$\max(Q(v) = |V(G')|)$}
            {
                $S \gets \pref{G'} \cup \{v\}$ for some $v$ that maximizes $Q(v)$\\
                \Return $S$
            }
            $i \gets \max(2, \varphi-|\pref{G'}|)$\\
            PDS $\gets$ \textbf{false}\\
            \While{\textbf{not} PDS}
            {
                \For{$C \subseteq A$ where $|C| = i$}
                {
                    $S \gets \pref{G'} \cup C$\\
                    \If{$S \cap R \neq \emptyset$ for each $R \in \jeff{G'}$}
                    {
                        \If{$\observed{G'}{S} = V(G')$}
                        {
                            PDS $\gets$ \textbf{true}\\
                            \Return $S$
                        }
                    }
                }
                $i \gets i+1$\\
            }
            \caption{The PDT algorithm for determining a minimum PDS of a \emph{connected} graph}
            \label{alg:connectedalgorithm}
        \end{algorithm}
        
\section{Run-time Analysis}\label{sec:survey}
    We now compare the JL-BW algorithm and the PDT in two ways: empirically, and with run-time examples on both random graphs with varying size and common IEEE test systems.

    \subsection{Empirical Comparison}
        For a graph $G$, we can compare the number of subsets strictly smaller than $\pdn{G}$ checked by both the JL-BW algorithm and the PDT.
        Let $N$ and $N'$ represent this number for the JL-BW algorithm and the PDT respectively:
        \begin{align*}
            N&:=\displaystyle\sum_{i=1}^{\pdn{G}-1}{\binom{|V(G)|}{i}}\\
            N'&:=1+a+\displaystyle\sum_{i=\max(2,\varphi)}^{\pdn{G}-p-1}\binom{a-\varphi}{i}2^{\varphi - p}
        \end{align*}
        where $G'$ is the contracted graph, $p=|\pref{G'}|$, $\varphi$ is as in Section~\ref{Section:guaranteed}, and $a$ is the number of active vertices with respect to the restricted power domination problem on $G'$ subject to $\pref{G'}$. 
        There is 1 case from checking if $\pref{G'}$ is a PDS and $a$ cases from determining if any vertex satisfies $Q(v)=|V(G')|$.
        We then check sets containing $i$ additional vertices from the set of active vertices that intersect non-trivially with each of the paired entrance vertices.
        Observe that $N'$ is often lower than $N$ due to the prevalence of preferred vertices, paired entrance vertices, and non-active vertices.
        
    \subsection{IEEE Test Systems}
\begin{figure}
            \centering
            \begin{subfigure}{0.49\textwidth}
                \resizebox{\textwidth}{!}{
                    \begin{tikzpicture}
                        \node [style=empty circle] (0) at (10, 6) {};
                        \node [style=empty circle] (1) at (10, 5) {};
                        \node [style=empty circle] (2) at (9, 5) {};
                        \node [style=empty circle] (3) at (9, 6) {};
                        \node [style=empty circle] (4) at (8, 6) {};
                        \node [style=empty circle] (5) at (8, 5) {};
                        \node [style=solid star] (6) at (7, 5) {};
                        \node [style=solid star] (7) at (8, 4) {};
                        \node [style=empty circle] (9) at (9, 3) {};
                        \node [style=empty circle] (11) at (7, 3) {};
                        \node [style=solid star] (12) at (7, 1) {};
                        \node [style=empty circle] (13) at (7, 2) {};
                        \node [style=empty circle] (14) at (8, 1) {};
                        \node [style=empty circle] (15) at (9, 1) {};
                        \node [style=empty circle] (16) at (8, 0) {};
                        \node [style=empty circle] (17) at (6, 5) {};
                        \node [style=empty circle] (18) at (6, 3) {};
                        \node [style=empty circle] (19) at (6, 0) {};
                        \node [style=empty circle] (20) at (5, 0) {};
                        \node [style=partial circle] (21) at (6, 1) {};
                        \node [style=partial circle] (22) at (5, 1) {};
                        \node [style=partial circle] (23) at (5, 3) {};
                        \node [style=partial circle] (24) at (4, 4) {};
                        \node [style=partial circle] (25) at (5, 5) {};
                        \node [style=partial circle] (26) at (3, 3) {};
                        \node [style=partial circle] (27) at (1, 3) {};
                        \node [style=partial circle] (28) at (3, 5) {};
                        \node [style=partial circle] (29) at (2, 5) {};
                        \node [style=empty circle] (30) at (2, 6) {};
                        \node [style=empty circle] (31) at (2, 4) {};
                        \node [style=solid circle] (32) at (1, 5) {};
                        \node [style=solid circle] (33) at (2, 2) {};
                        \node [style=partial diamond] (34) at (10, 4) {};
                        \node [style=partial diamond] (35) at (2, 1) {};
                        \node [style=partial diamond] (36) at (3, 1) {};
                        \node [style=partial diamond] (37) at (4, 1) {};
                        \node [style=empty circle] (38) at (1, 2) {};
                        \node [style=empty circle] (39) at (0, 3) {};
                        \node [style=partial diamond] (40) at (9, 4) {};
                        \draw [style=solid edge] (4) to (3);
                        \draw [style=solid edge] (3) to (0);
                        \draw [style=solid edge] (3) to (2);
                        \draw [style=solid edge] (2) to (1);
                        \draw [style=solid edge] (5) to (2);
                        \draw [style=solid edge] (6) to (4);
                        \draw [style=solid edge] (6) to (5);
                        \draw [style=solid edge] (7) to (9);
                        \draw [style=solid edge] (7) to (6);
                        \draw [style=solid edge] (6) to (11);
                        \draw [style=solid edge] (15) to (14);
                        \draw [style=solid edge] (14) to (16);
                        \draw [style=solid edge] (14) to (12);
                        \draw [style=solid edge] (12) to (13);
                        \draw [style=solid edge] (13) to (11);
                        \draw [style=solid edge] (17) to (6);
                        \draw [style=solid edge] (17) to (25);
                        \draw [style=solid edge] (25) to (24);
                        \draw [style=solid edge] (24) to (23);
                        \draw [style=solid edge] (23) to (18);
                        \draw [style=solid edge] (11) to (18);
                        \draw [style=solid edge] (21) to (19);
                        \draw [style=solid edge] (20) to (22);
                        \draw [style=solid edge] (22) to (21);
                        \draw [style=solid edge] (21) to (12);
                        \draw [style=solid edge] (22) to (23);
                        \draw [style=solid edge] (27) to (26);
                        \draw [style=solid edge] (26) to (24);
                        \draw [style=solid edge] (28) to (25);
                        \draw [style=solid edge] (29) to (28);
                        \draw [style=solid edge] (31) to (28);
                        \draw [style=solid edge] (29) to (30);
                        \draw [style=solid edge] (32) to (29);
                        \draw [style=solid edge] (31) to (32);
                        \draw [style=solid edge] (32) to (27);
                        \draw [style=solid edge] (22) to (37);
                        \draw [style=solid edge] (37) to (36);
                        \draw [style=solid edge] (36) to (35);
                        \draw [style=solid edge] (35) to (33);
                        \draw [style=solid edge] (33) to (26);
                        \draw [style=solid edge] (39) to (27);
                        \draw [style=solid edge] (27) to (38);
                        \draw [style=solid edge] (38) to (33);
                        \draw [style=solid edge] (7) to (40);
                        \draw [style=solid edge] (40) to (34);
                    \end{tikzpicture}
                }
                \caption{New England Power Grid as a graph.}
                \label{fig:newengland}
            \end{subfigure}
            \hfill
            \begin{subfigure}{0.49\textwidth}
                \resizebox{\textwidth}{!}{
                    \begin{tikzpicture}
                        \node [style=empty circle] (0) at (10, 6) {};
                        \node [style=empty circle] (1) at (10, 5) {};
                        \node [style=empty circle] (2) at (9, 5) {};
                        \node [style=empty circle] (3) at (9, 6) {};
                        \node [style=empty circle] (4) at (8, 6) {};
                        \node [style=empty circle] (5) at (8, 5) {};
                        \node [style=solid star] (6) at (7, 5) {};
                        \node [style=solid star] (7) at (8, 4) {};
                        \node [style=empty circle] (8) at (9, 4) {};
                        \node [style=empty circle] (9) at (9, 3) {};
                        \node [style=empty circle] (11) at (7, 3) {};
                        \node [style=solid star] (12) at (7, 1) {};
                        \node [style=empty circle] (13) at (7, 2) {};
                        \node [style=empty circle] (14) at (8, 1) {};
                        \node [style=empty circle] (15) at (9, 1) {};
                        \node [style=empty circle] (16) at (8, 0) {};
                        \node [style=empty circle] (17) at (6, 5) {};
                        \node [style=empty circle] (18) at (6, 3) {};
                        \node [style=empty circle] (19) at (6, 0) {};
                        \node [style=empty circle] (20) at (5, 0) {};
                        \node [style=partial circle] (21) at (6, 1) {};
                        \node [style=partial circle] (22) at (5, 1) {};
                        \node [style=partial circle] (23) at (5, 3) {};
                        \node [style=partial circle] (24) at (4, 4) {};
                        \node [style=partial circle] (25) at (5, 5) {};
                        \node [style=partial circle] (26) at (3, 3) {};
                        \node [style=partial circle] (27) at (1, 3) {};
                        \node [style=partial circle] (28) at (3, 5) {};
                        \node [style=partial circle] (29) at (2, 5) {};
                        \node [style=empty circle] (30) at (2, 6) {};
                        \node [style=empty circle] (31) at (2, 4) {};
                        \node [style=solid circle] (32) at (1, 5) {};
                        \node [style=solid circle] (33) at (2, 2) {};
                        \node [style=empty circle] (38) at (1, 2) {};
                        \node [style=empty circle] (39) at (0, 3) {};
                        \node [style=empty circle] (40) at (3.25, 1.25) {};
                        \draw [style=solid edge] (4) to (3);
                        \draw [style=solid edge] (3) to (0);
                        \draw [style=solid edge] (3) to (2);
                        \draw [style=solid edge] (2) to (1);
                        \draw [style=solid edge] (5) to (2);
                        \draw [style=solid edge] (6) to (4);
                        \draw [style=solid edge] (6) to (5);
                        \draw [style=solid edge] (8) to (7);
                        \draw [style=solid edge] (7) to (9);
                        \draw [style=solid edge] (7) to (6);
                        \draw [style=solid edge] (6) to (11);
                        \draw [style=solid edge] (15) to (14);
                        \draw [style=solid edge] (14) to (16);
                        \draw [style=solid edge] (14) to (12);
                        \draw [style=solid edge] (12) to (13);
                        \draw [style=solid edge] (13) to (11);
                        \draw [style=solid edge] (17) to (6);
                        \draw [style=solid edge] (17) to (25);
                        \draw [style=solid edge] (25) to (24);
                        \draw [style=solid edge] (24) to (23);
                        \draw [style=solid edge] (23) to (18);
                        \draw [style=solid edge] (11) to (18);
                        \draw [style=solid edge] (21) to (19);
                        \draw [style=solid edge] (20) to (22);
                        \draw [style=solid edge] (22) to (21);
                        \draw [style=solid edge] (21) to (12);
                        \draw [style=solid edge] (22) to (23);
                        \draw [style=solid edge] (27) to (26);
                        \draw [style=solid edge] (26) to (24);
                        \draw [style=solid edge] (28) to (25);
                        \draw [style=solid edge] (29) to (28);
                        \draw [style=solid edge] (31) to (28);
                        \draw [style=solid edge] (29) to (30);
                        \draw [style=solid edge] (32) to (29);
                        \draw [style=solid edge] (31) to (32);
                        \draw [style=solid edge] (32) to (27);
                        \draw [style=solid edge] (33) to (26);
                        \draw [style=solid edge] (39) to (27);
                        \draw [style=solid edge] (27) to (38);
                        \draw [style=solid edge] (38) to (33);
                        \draw [style=solid edge] (40) to (33);
                        \draw [style=solid edge] (40) to (22);
                    \end{tikzpicture}
                }
                \caption{Contracted New England Power Grid as a graph.}
                \label{fig:contractednewengland}
            \end{subfigure}
            \caption{New England power grid as a graph with 39 vertices in (a). The PDT removes the dark grey diamond vertices to contract the graph into the graph in (b). Further, the PDT locates 3 preferred vertices (black stars) and 11 active vertices 3 (light grey and black circles). The minimum PDS of size 5 returned by the PDT is indicated by solid black vertices.}
            \label{Figure:NewEngland}
        \end{figure}
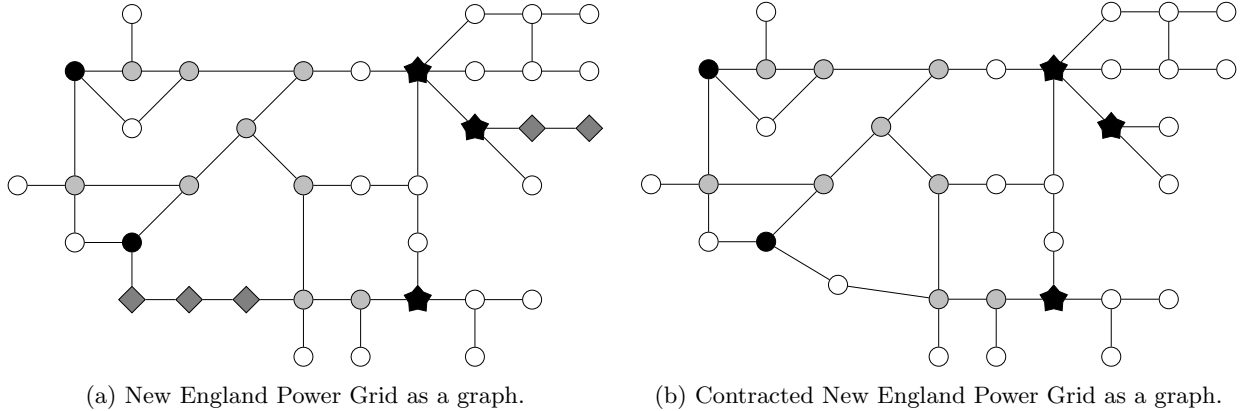

        Returning to the original problem of the 2003 power grid failure, we consider the IEEE 39 bus test system that represents a historic model of the New England power grid as available in the pandapower Python module \cite{Thurner2018}.
        This graph is displayed in Figure~\ref{fig:newengland} and has power domination number 5.
        The JL-BW algorithm evaluates the $N=$~92,170 subsets to determine $\pdn{G} > 4$ and locates a minimum PDS of size 5 in an average of 3.566 seconds.
        The PDT contracts to a graph on 36 vertices and finds 3 preferred vertices.
        The contracted graph is shown in Figure~\ref{fig:contractednewengland}.
        By considering the restricted power domination problem on $G'$ subject to $\pref{G'}$, the PDT indicates 11 active vertices. 
        The PDT evaluates $N'=12$ subsets to determine $\pdn{G} > 4$ and locates a minimum PDS of size 5 in an average of $2.673 \times 10 ^ {-3}$ seconds.
        
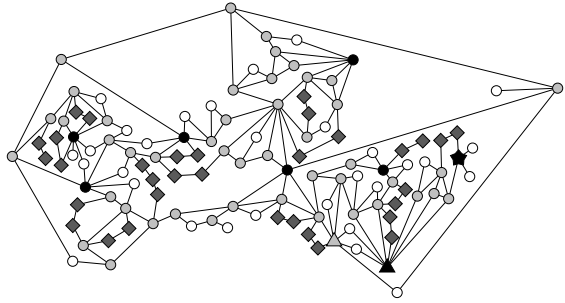
\begin{figure}
            \centering
            \resizebox{0.45\textwidth}{!}{  
            \begin{tikzpicture}
      \draw
        (-13.856, 3.383) node[shape=diamond, draw=black, fill=black, fill opacity=.65, text opacity=0] (1){00}
        (-13.508, 1.353) node[shape=diamond, draw=black, fill=black, fill opacity=.65, text opacity=0] (2){00}
        (-13.323, 4.611) node[shape=circle, draw=black, fill=black, fill opacity=.25, text opacity=0] (3){00}
        (-10.575, 6.269) node[shape=circle, draw=black, fill=black, fill opacity=0, text opacity=0] (4){00}
        (-12.563, 6.794) node[shape=circle, draw=black, fill=black, fill opacity=.25, text opacity=0] (5){00}
        (-12.416, 5.271) node[shape=diamond, draw=black, fill=black, fill opacity=.65, text opacity=0] (6){00}
        (-11.403, 4.81) node[shape=diamond, draw=black, fill=black, fill opacity=.65, text opacity=0] (7){00}
        (-14.278, 4.792) node[shape=circle, draw=black, fill=black, fill opacity=.25, text opacity=0] (8){00}
        (-15.148, 2.978) node[shape=diamond, draw=black, fill=black, fill opacity=.65, text opacity=0] (9){00}
        (-14.663, 1.853) node[shape=diamond, draw=black, fill=black, fill opacity=.65, text opacity=0] (10){00}
        (-10.116, 4.641) node[shape=circle, draw=black, fill=black, fill opacity=.25, text opacity=0] (11){00}
        (-12.597, 3.444) node[shape=circle, draw=black, fill=black, fill opacity=1, text opacity=0] (12){00}
        (-8.671, 3.919) node[shape=circle, draw=black, fill=black, fill opacity=0, text opacity=0] (13){00}
        (-11.362, 2.414) node[shape=circle, draw=black, fill=black, fill opacity=0, text opacity=0] (14){00}
        (-9.969, 3.222) node[shape=circle, draw=black, fill=black, fill opacity=.25, text opacity=0] (15){00}
        (-11.841, 1.453) node[shape=circle, draw=black, fill=black, fill opacity=0, text opacity=0] (16){00}
        (-11.723, -0.263) node[shape=circle, draw=black, fill=black, fill opacity=1, text opacity=0] (17){00}
        (-8.952, 0.828) node[shape=circle, draw=black, fill=black, fill opacity=0, text opacity=0] (18){00}
        (-8.393, 2.289) node[shape=circle, draw=black, fill=black, fill opacity=.25, text opacity=0] (19){00}
        (-7.526, 1.352) node[shape=diamond, draw=black, fill=black, fill opacity=.65, text opacity=0] (20){00}
        (-6.756, 0.296) node[shape=diamond, draw=black, fill=black, fill opacity=.65, text opacity=0] (21){00}
        (-6.403, -0.809) node[shape=diamond, draw=black, fill=black, fill opacity=.65, text opacity=0] (22){00}
        (-6.733, -2.955) node[shape=circle, draw=black, fill=black, fill opacity=.25, text opacity=0] (23){00}
        (-5.085, -2.244) node[shape=circle, draw=black, fill=black, fill opacity=.25, text opacity=0] (24){00}
        (-9.862, -5.992) node[shape=circle, draw=black, fill=black, fill opacity=.25, text opacity=0] (25){00}
        (-12.661, -5.732) node[shape=circle, draw=black, fill=black, fill opacity=0, text opacity=0] (26){00}
        (-11.88, -4.569) node[shape=circle, draw=black, fill=black, fill opacity=.25, text opacity=0] (27){00}
        (-12.652, -3.093) node[shape=diamond, draw=black, fill=black, fill opacity=.65, text opacity=0] (28){00}
        (-12.309, -1.563) node[shape=diamond, draw=black, fill=black, fill opacity=.65, text opacity=0] (29){00}
        (-17.113, 2.006) node[shape=circle, draw=black, fill=black, fill opacity=.25, text opacity=0] (30){00}
        (-9.859, -0.971) node[shape=circle, draw=black, fill=black, fill opacity=.25, text opacity=0] (31){00}
        (-8.751, -1.83) node[shape=circle, draw=black, fill=black, fill opacity=.25, text opacity=0] (32){00}
        (-7.187, 2.947) node[shape=circle, draw=black, fill=black, fill opacity=0, text opacity=0] (33){00}
        (-6.554, 1.937) node[shape=circle, draw=black, fill=black, fill opacity=.25, text opacity=0] (34){00}
        (-3.437, 2.096) node[shape=diamond, draw=black, fill=black, fill opacity=.65, text opacity=0] (35){00}
        (-4.968, 1.982) node[shape=diamond, draw=black, fill=black, fill opacity=.65, text opacity=0] (36){00}
        (-4.452, 3.407) node[shape=circle, draw=black, fill=black, fill opacity=1, text opacity=0] (37){00}
        (-13.505, 9.162) node[shape=circle, draw=black, fill=black, fill opacity=.25, text opacity=0] (38){00}
        (-4.386, 4.955) node[shape=circle, draw=black, fill=black, fill opacity=0, text opacity=0] (39){00}
        (-2.429, 3.106) node[shape=circle, draw=black, fill=black, fill opacity=.25, text opacity=0] (40){00}
        (-2.448, 5.737) node[shape=circle, draw=black, fill=black, fill opacity=0, text opacity=0] (41){00}
        (-1.359, 4.704) node[shape=circle, draw=black, fill=black, fill opacity=.25, text opacity=0] (42){00}
        (-5.15, 0.577) node[shape=diamond, draw=black, fill=black, fill opacity=.65, text opacity=0] (43){00}
        (-3.12, 0.692) node[shape=diamond, draw=black, fill=black, fill opacity=.65, text opacity=0] (44){00}
        (-1.498, 2.416) node[shape=circle, draw=black, fill=black, fill opacity=.25, text opacity=0] (45){00}
        (-0.3, 1.513) node[shape=circle, draw=black, fill=black, fill opacity=.25, text opacity=0] (46){00}
        (1.711, 2.072) node[shape=circle, draw=black, fill=black, fill opacity=.25, text opacity=0] (47){00}
        (0.881, 3.421) node[shape=circle, draw=black, fill=black, fill opacity=0, text opacity=0] (48){00}
        (2.479, 5.852) node[shape=circle, draw=black, fill=black, fill opacity=.25, text opacity=0] (49){00}
        (4.064, 2.004) node[shape=diamond, draw=black, fill=black, fill opacity=.65, text opacity=0] (50){00}
        (4.615, 3.428) node[shape=circle, draw=black, fill=black, fill opacity=.25, text opacity=0] (51){00}
        (4.751, 5.162) node[shape=diamond, draw=black, fill=black, fill opacity=.65, text opacity=0] (52){00}
        (4.4, 6.428) node[shape=diamond, draw=black, fill=black, fill opacity=.65, text opacity=0] (53){00}
        (4.64, 7.844) node[shape=circle, draw=black, fill=black, fill opacity=.25, text opacity=0] (54){00}
        (6.449, 7.606) node[shape=circle, draw=black, fill=black, fill opacity=.25, text opacity=0] (55){00}
        (6.853, 5.833) node[shape=circle, draw=black, fill=black, fill opacity=.25, text opacity=0] (56){00}
        (6.939, 3.466) node[shape=diamond, draw=black, fill=black, fill opacity=.65, text opacity=0] (57){00}
        (5.96, 4.183) node[shape=circle, draw=black, fill=black, fill opacity=0, text opacity=0] (58){00}
        (8.018, 9.131) node[shape=circle, draw=black, fill=black, fill opacity=1, text opacity=0] (59){00}
        (3.658, 8.708) node[shape=circle, draw=black, fill=black, fill opacity=.25, text opacity=0] (60){00}
        (2.302, 9.739) node[shape=circle, draw=black, fill=black, fill opacity=.25, text opacity=0] (61){00}
        (2.024, 7.751) node[shape=circle, draw=black, fill=black, fill opacity=.25, text opacity=0] (62){00}
        (3.865, 10.598) node[shape=circle, draw=black, fill=black, fill opacity=0, text opacity=0] (63){00}
        (1.621, 10.858) node[shape=circle, draw=black, fill=black, fill opacity=.25, text opacity=0] (64){00}
        (-1.006, 12.956) node[shape=circle, draw=black, fill=black, fill opacity=.25, text opacity=0] (65){00}
        (-0.825, 6.904) node[shape=circle, draw=black, fill=black, fill opacity=.25, text opacity=0] (66){00}
        (0.66, 8.434) node[shape=circle, draw=black, fill=black, fill opacity=0, text opacity=0] (67){00}
        (23.104, 7.043) node[shape=circle, draw=black, fill=black, fill opacity=.25, text opacity=0] (68){00}
        (3.166, 1.003) node[shape=circle, draw=black, fill=black, fill opacity=1, text opacity=0] (69){00}
        (-0.833, -1.684) node[shape=circle, draw=black, fill=black, fill opacity=.25, text opacity=0] (70){00}
        (-2.384, -2.719) node[shape=circle, draw=black, fill=black, fill opacity=.25, text opacity=0] (71){00}
        (-3.9, -3.137) node[shape=circle, draw=black, fill=black, fill opacity=0, text opacity=0] (72){00}
        (-1.24, -3.257) node[shape=circle, draw=black, fill=black, fill opacity=0, text opacity=0] (73){00}
        (0.84, -2.342) node[shape=circle, draw=black, fill=black, fill opacity=0, text opacity=0] (74){00}
        (2.756, -1.163) node[shape=circle, draw=black, fill=black, fill opacity=.25, text opacity=0] (75){00}
        (3.626, -2.806) node[shape=diamond, draw=black, fill=black, fill opacity=.65, text opacity=0] (76){00}
        (5.507, -2.468) node[shape=circle, draw=black, fill=black, fill opacity=.25, text opacity=0] (77){00}
        (4.719, -3.752) node[shape=diamond, draw=black, fill=black, fill opacity=.65, text opacity=0] (78){00}
        (5.423, -4.751) node[shape=diamond, draw=black, fill=black, fill opacity=.65, text opacity=0] (79){00}
        (6.619, -4.27) node[shape=regular polygon, regular polygon sides=3, draw=black, fill=black, fill opacity=.25, text opacity=0] (80){0}
        (11.266, -8.04) node[shape=circle, draw=black, fill=black, fill opacity=0, text opacity=0] (81){00}
        (5.012, 0.574) node[shape=circle, draw=black, fill=black, fill opacity=.25, text opacity=0] (82){00}
        (7.851, 1.429) node[shape=circle, draw=black, fill=black, fill opacity=.25, text opacity=0] (83){00}
        (9.454, 2.293) node[shape=circle, draw=black, fill=black, fill opacity=0, text opacity=0] (84){00}
        (10.243, 0.993) node[shape=circle, draw=black, fill=black, fill opacity=1, text opacity=0] (85){00}
        (11.623, 2.435) node[shape=diamond, draw=black, fill=black, fill opacity=.65, text opacity=0] (86){00}
        (13.156, 3.149) node[shape=diamond, draw=black, fill=black, fill opacity=.65, text opacity=0] (87){00}
        (12.03, 1.391) node[shape=circle, draw=black, fill=black, fill opacity=0, text opacity=0] (88){00}
        (10.932, 0.116) node[shape=circle, draw=black, fill=black, fill opacity=.25, text opacity=0] (89){00}
        (11.871, -0.445) node[shape=diamond, draw=black, fill=black, fill opacity=.65, text opacity=0] (90){00}
        (11.128, -1.486) node[shape=diamond, draw=black, fill=black, fill opacity=.65, text opacity=0] (91){00}
        (9.835, -1.546) node[shape=circle, draw=black, fill=black, fill opacity=.25, text opacity=0] (92){00}
        (9.799, -0.244) node[shape=circle, draw=black, fill=black, fill opacity=0, text opacity=0] (93){00}
        (8.038, -2.265) node[shape=circle, draw=black, fill=black, fill opacity=.25, text opacity=0] (94){00}
        (8.349, 0.564) node[shape=circle, draw=black, fill=black, fill opacity=0, text opacity=0] (95){00}
        (7.118, 0.36) node[shape=circle, draw=black, fill=black, fill opacity=.25, text opacity=0] (96){00}
        (6.108, -1.543) node[shape=circle, draw=black, fill=black, fill opacity=0, text opacity=0] (97){00}
        (8.261, -4.84) node[shape=circle, draw=black, fill=black, fill opacity=0, text opacity=0] (98){00}
        (7.766, -3.328) node[shape=circle, draw=black, fill=black, fill opacity=0, text opacity=0] (99){00}
        (10.532, -6.2) node[shape=regular polygon, regular polygon sides=3, draw=black, fill=black, fill opacity=1, text opacity=0] (100){0}
        (10.874, -3.974) node[shape=diamond, draw=black, fill=black, fill opacity=0.65, text opacity=0] (101){00}
        (10.741, -2.462) node[shape=diamond, draw=black, fill=black, fill opacity=0.65, text opacity=0] (102){00}
        (15.065, -1.115) node[shape=circle, draw=black, fill=black, fill opacity=.25, text opacity=0] (103){00}
        (13.987, -0.736) node[shape=circle, draw=black, fill=black, fill opacity=.25, text opacity=0] (104){00}
        (14.545, 0.815) node[shape=circle, draw=black, fill=black, fill opacity=.25, text opacity=0] (105){00}
        (12.758, -0.923) node[shape=circle, draw=black, fill=black, fill opacity=.25, text opacity=0] (106){00}
        (13.306, 1.599) node[shape=circle, draw=black, fill=black, fill opacity=0, text opacity=0] (107){00}
        (14.481, 3.187) node[shape=diamond, draw=black, fill=black, fill opacity=0.65, text opacity=0] (108){00}
        (15.694, 3.749) node[shape=diamond, draw=black, fill=black, fill opacity=0.65, text opacity=0] (109){00}
        (15.794, 1.876) node[shape=star, draw=black, fill=black, fill opacity=1, text opacity=0] (110){00}
        (16.617, 0.517) node[shape=circle, draw=black, fill=black, fill opacity=0, text opacity=0] (111){00}
        (16.83, 2.622) node[shape=circle, draw=black, fill=black, fill opacity=0, text opacity=0] (112){00}
        (-8.119, -0.008) node[shape=circle, draw=black, fill=black, fill opacity=0, text opacity=0] (113){00}
        (-8.514, -3.31) node[shape=diamond, draw=black, fill=black, fill opacity=0.65, text opacity=0] (114){00}
        (-10.026, -4.21) node[shape=diamond, draw=black, fill=black, fill opacity=0.65, text opacity=0] (115){00}
        (18.585, 6.852) node[shape=circle, draw=black, fill=black, fill opacity=0, text opacity=0] (116){00}
        (-12.654, 2.094) node[shape=circle, draw=black, fill=black, fill opacity=0, text opacity=0] (117){00}
        (2.46, -2.505) node[shape=diamond, draw=black, fill=black, fill opacity=0.65, text opacity=0] (118){00};
      \begin{scope}[-]
        \draw (1) to (2);
        \draw (1) to (3);
        \draw (2) to (12);
        \draw (3) to (5);
        \draw (3) to (12);
        \draw (4) to (5);
        \draw (4) to (11);
        \draw (5) to (6);
        \draw (5) to (11);
        \draw (5) to (8);
        \draw (6) to (7);
        \draw (7) to (12);
        \draw (8) to (9);
        \draw (8) to (30);
        \draw (9) to (10);
        \draw (11) to (12);
        \draw (11) to (13);
        \draw (12) to (14);
        \draw (12) to (16);
        \draw (12) to (117);
        \draw (13) to (15);
        \draw (14) to (15);
        \draw (15) to (17);
        \draw (15) to (19);
        \draw (15) to (33);
        \draw (16) to (17);
        \draw (17) to (18);
        \draw (17) to (31);
        \draw (17) to (113);
        \draw (17) to (30);
        \draw (18) to (19);
        \draw (19) to (20);
        \draw (19) to (34);
        \draw (20) to (21);
        \draw (21) to (22);
        \draw (22) to (23);
        \draw (23) to (24);
        \draw (23) to (25);
        \draw (23) to (32);
        \draw (24) to (70);
        \draw (24) to (72);
        \draw (25) to (27);
        \draw (25) to (26);
        \draw (26) to (30);
        \draw (27) to (28);
        \draw (27) to (32);
        \draw (27) to (115);
        \draw (28) to (29);
        \draw (29) to (31);
        \draw (30) to (38);
        \draw (31) to (32);
        \draw (32) to (113);
        \draw (32) to (114);
        \draw (33) to (37);
        \draw (34) to (36);
        \draw (34) to (37);
        \draw (34) to (43);
        \draw (35) to (36);
        \draw (35) to (37);
        \draw (37) to (39);
        \draw (37) to (40);
        \draw (37) to (38);
        \draw (38) to (65);
        \draw (39) to (40);
        \draw (40) to (41);
        \draw (40) to (42);
        \draw (41) to (42);
        \draw (42) to (49);
        \draw (43) to (44);
        \draw (44) to (45);
        \draw (45) to (46);
        \draw (45) to (49);
        \draw (46) to (47);
        \draw (46) to (48);
        \draw (47) to (49);
        \draw (47) to (69);
        \draw (48) to (49);
        \draw (49) to (50);
        \draw (49) to (51);
        \draw (49) to (54);
        \draw (49) to (66);
        \draw (49) to (69);
        \draw (50) to (57);
        \draw (51) to (52);
        \draw (51) to (58);
        \draw (52) to (53);
        \draw (53) to (54);
        \draw (54) to (55);
        \draw (54) to (56);
        \draw (54) to (59);
        \draw (55) to (56);
        \draw (55) to (59);
        \draw (56) to (57);
        \draw (56) to (58);
        \draw (56) to (59);
        \draw (59) to (60);
        \draw (59) to (61);
        \draw (59) to (63);
        \draw (60) to (61);
        \draw (60) to (62);
        \draw (61) to (62);
        \draw (61) to (64);
        \draw (62) to (66);
        \draw (62) to (67);
        \draw (63) to (64);
        \draw (64) to (65);
        \draw (65) to (66);
        \draw (65) to (68);
        \draw (66) to (67);
        \draw (68) to (69);
        \draw (68) to (81);
        \draw (68) to (116);
        \draw (69) to (70);
        \draw (69) to (75);
        \draw (69) to (77);
        \draw (70) to (71);
        \draw (70) to (74);
        \draw (70) to (75);
        \draw (71) to (72);
        \draw (71) to (73);
        \draw (74) to (75);
        \draw (75) to (77);
        \draw (75) to (118);
        \draw (76) to (77);
        \draw (76) to (118);
        \draw (77) to (78);
        \draw (77) to (80);
        \draw (77) to (82);
        \draw (78) to (79);
        \draw (79) to (80);
        \draw (80) to (96);
        \draw (80) to (97);
        \draw (80) to (98);
        \draw (80) to (99);
        \draw (80) to (81);
        \draw (82) to (83);
        \draw (82) to (96);
        \draw (83) to (84);
        \draw (83) to (85);
        \draw (84) to (85);
        \draw (85) to (86);
        \draw (85) to (88);
        \draw (85) to (89);
        \draw (86) to (87);
        \draw (88) to (89);
        \draw (89) to (90);
        \draw (89) to (92);
        \draw (90) to (91);
        \draw (91) to (92);
        \draw (92) to (93);
        \draw (92) to (94);
        \draw (92) to (100);
        \draw (92) to (102);
        \draw (93) to (94);
        \draw (94) to (95);
        \draw (94) to (96);
        \draw (94) to (100);
        \draw (95) to (96);
        \draw (96) to (97);
        \draw (98) to (100);
        \draw (99) to (100);
        \draw (100) to (101);
        \draw (100) to (103);
        \draw (100) to (104);
        \draw (100) to (106);
        \draw (101) to (102);
        \draw (103) to (104);
        \draw (103) to (105);
        \draw (103) to (110);
        \draw (104) to (105);
        \draw (105) to (106);
        \draw (105) to (107);
        \draw (105) to (108);
        \draw (106) to (107);
        \draw (108) to (109);
        \draw (109) to (110);
        \draw (110) to (111);
        \draw (110) to (112);
        \draw (114) to (115);
      \end{scope}
    \end{tikzpicture}
            }
            \caption{IEEE 118 vertex test system as a graph. The PDT contracts the dark grey diamond vertices to yield $G'$ with 101 vertices. The PDT indicates 1 preferred vertex (black star), 1 set of paired entrance vertices (triangles), and 54 active vertices (light grey and black circles and triangles). The minimum PDS of size 8 returned by the PDT is indicated by solid black vertices.}
            \label{fig:118}
        \end{figure}

        The IEEE 118 bus test system \cite{IEEE118} shown in Figure~\ref{fig:118} is also available in the pandapower module and has $\pdn{G}=8$.
        The JL-BW algorithm would need to evaluate the $N\approx 5.620\times10^{10}$ subsets to determine $\pdn{G}$, however the search was terminated after a week and no minimum PDS was located.
        The PDT locates the single preferred vertex, contracts the graph to a graph with 115 vertices, locates 54 active vertices, and finds one set of paired entrance vertices.
        By leveraging the PDT and 32 threads, $N'\approx4.650\times10^7$ subsets are evaluated to determine $\pdn{G} > 7$ and to locate a minimum PDS of size 8 in $1.098\times10^2$ seconds.
        This exemplifies the drastic run-time improvements of the PDT over the JL-BW algorithm.

        Table~\ref{tab:ieeeruntimecomparison} shows the average run-time on various other IEEE test systems up to 300 vertices.
        In the appendix, Table~\ref{tab:variousieeegraphs} details $|V(G)|$, $|V(G')|$, $|\pref{G'}|$, and $|\jeff{G'}|$ for other IEEE test systems up to the 1354 vertex system.
        \begin{table}[ht]
            \centering
            \begin{tabular}{|c|c||l|l|}
                \hline
                \begin{tabular}[c]{@{}l@{}}number of\\ vertices\end{tabular} & $\pdn{G}$ & \begin{tabular}[c]{@{}l@{}}average JL-BW\\ algorithm time\end{tabular} & average PDT time         \\ \hline
                5                       & 1                                       & $\mathbf{2.087 \times 10 ^ {-5}}$ s                                              & $1.526 \times 10 ^ {-4}$ s \\ \hline
                6                       & 1                                     & $\mathbf{2.791 \times 10 ^ {-5}}$ s                                              & $1.643 \times 10 ^ {-4}$ s \\ \hline
                9                       & 1                                    & $\mathbf{9.541 \times 10 ^ {-5}}$ s                                              & $7.172 \times 10 ^ {-4}$ s \\ \hline
                11                      & 2                                    & $\mathbf{1.518 \times 10 ^ {-4}}$ s                                              & $3.194 \times 10 ^ {-4}$ s \\ \hline
                14                      & 2                                    & $\mathbf{1.806 \times 10 ^ {-4}}$ s                                              & $1.023 \times 10 ^ {-3}$ s \\ \hline
                24                      & 3                                    & $1.161 \times 10 ^ {-2}$ s                                              & $\mathbf{2.043 \times 10 ^ {-3}}$ s \\ \hline
                30                      & 3                                    & $8.830 \times 10 ^ {-3}$ s                                              & $\mathbf{1.920 \times 10 ^ {-3}}$ s \\ \hline
                30                      & 3                                    & $9.259 \times 10 ^ {-3}$ s                                              & $\mathbf{1.954 \times 10 ^ {-3}}$ s \\ \hline
                33                      & 1                                    & $5.945 \times 10 ^ {-4}$ s                                              & $\mathbf{5.368 \times 10 ^ {-4}}$ s \\ \hline
                39                      & 5                                   & $3.567$ s                                                               & $\mathbf{2.733 \times 10 ^ {-3}}$ s \\ \hline
                57                      & 3                                    & $2.352 \times 10 ^ {-1}$ s                                              & $\mathbf{1.064 \times 10 ^ {-2}}$ s \\ \hline
                89                      &  5                                   & $2.023 \times 10 ^ {2}$ s                                               & $\mathbf{7.776 \times 10 ^ {-3}}$ s \\ \hline
                118                     &  8                                   & $>1$ week                                                         & $\mathbf{1.098 \times 10 ^ {2}}$ s  \\ \hline
                145                     &  13                                   & N/A                                                         & $\mathbf{6.495 \times 10 ^ {5}}$ s$^*$  \\ \hline
                200                     &  20                                   & N/A                                                         & $\mathbf{4.062 \times 10 ^ {-2}}$ s \\ \hline
                300                     &  ?                                   & N/A                                                         & $>1$ week            \\ \hline
            \end{tabular}
            \caption{Average run-time for determining $\pdn{G}$ for various IEEE test systems. The 145 vertex test system, due to it's long run-time, was only tested once where each other graph was tested 20 times unless terminated early.} 
            \label{tab:ieeeruntimecomparison}
        \end{table}
        
    \subsection{On \ER Random Graphs With Varying Size}\label{sec:carlo}    
        We now compare the run-time of the PDT to the JL-BW algorithm on \ER random graphs.
        An \emph{\ER random graph} on $n$ vertices is a graph resulting from adding edges between each pair of distinct vertices with a predetermined probability \cite{Erdoes1959}.
        In this paper, we will consider \ER random graphs with edge probability of 0.05 and a varying number of vertices.
        While it is not guaranteed that an \ER random graph is connected, we generate \ER random graphs until the resulting graph is connected and store the connected graph.
        To collect the following run-time data, we determined $\pdn{G}$ for each graph 20 times and used the average run-time as the estimated run-time for each graph.
    
\begin{figure}[ht]
            \centering
            \includegraphics[width=\textwidth]{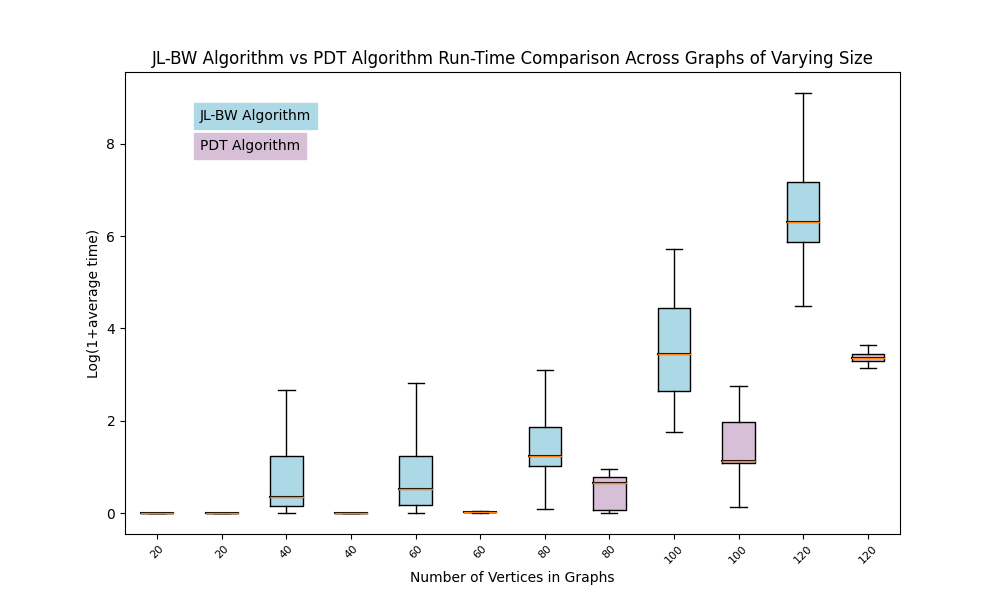}
            \caption{Boxplots of time to determine $\pdn{G}$ grouped by algorithm and $|V(G)|$ for each graph $G$ in the \ER data set including 100 graphs with 20, 40, 60, 80, 100, and 120 vertices each.}
            \label{fig:Order}
        \end{figure}

        We investigated the impact graph order has on time to find a minimum PDS by testing on connected \ER random graphs with 20, 40, 60, 80, 100, and 120 vertices.
        This data set is available in graph6 format \cite{Graph6} alongside the PDT and includes a total of 600 \ER random graphs (100 of each order).
        As expected, when the order of the graph increases, so does the time required to find a minimum PDS.
        This difference in run-time is shown in Figure~\ref{fig:Order} and Table~\ref{tab:summarystats} gives mean and median values.

        In addition to the faster run-times, the PDT yields less variance in run-time.
        \begin{table}[ht]
            \caption{Average run-time on \ER random graphs with varying number of vertices}
            \label{tab:summarystats}
            \centering
            \resizebox{0.95\textwidth}{!}{
                \begin{tabular}{|cl|l|l|l|l|l|l|}
                    \hline
                    \multicolumn{2}{|l|}{Algorithm} & $|V(G)| = 20$ & $|V(G)| = 40$ & $|V(G)| = 60$ & $|V(G)| = 80$ & $|V(G)| = 100$ & $|V(G)| = 120$ \\ \hline
                    \multicolumn{1}{|c|}{\multirow{2}{*}{mean}}   & JL-BW & $4.550 \times 10 ^ {-3}$ & $4.646$ & $6.361$ & $1.009 \times 10 ^ 1$ & $7.216 \times 10 ^ 1$ & $1.108 \times 10 ^ {2}$ \\ \cline{2-8} 
                    \multicolumn{1}{|c|}{}                        & PDT  & $\mathbf{6.692 \times 10 ^ {-4}}$ & $\mathbf{4.847 \times 10 ^ {-3}}$ & $\mathbf{7.233 \times 10 ^ {-2}}$ & $\mathbf{6.999 \times 10 ^ {-1}}$ & $\mathbf{4.077}$  & $\mathbf{3.935 \times 10 ^ {1}}$  \\ \hline
                    \multicolumn{1}{|c|}{\multirow{2}{*}{median}} & JL-BW & $2.661 \times 10 ^ {-3}$ & $4.218 \times 10 ^ {-1}$ & $7.004 \times 10 ^ {-1}$ & $2.416$ & $3.009 \times 10 ^ {1}$ & $5.483 \times 10 ^ {2}$ \\ \cline{2-8} 
                    \multicolumn{1}{|c|}{}                        & PDT  & $\mathbf{5.833 \times 10 ^ {-4}}$ & $\mathbf{2.959 \times 10 ^ {-3}}$  & $\mathbf{2.034 \times 10 ^ {-2}}$ & $\mathbf{9.241 \times 10 ^ {-1}}$ & $\mathbf{2.110}$  & $\mathbf{2.793 \times 10 ^{1}}$  \\ \hline
                \end{tabular}
            }
        \end{table}
        
\section{Concluding Remarks}\label{sec:concluding}

    \subsection{Using the Power Domination Toolbox}\label{sec:using}
        The primary interface functions of the PDT are the following:
        \begin{enumerate}
            \item
                \verb|PDT_pdn(Input_graph, Number_workers) -> int|
                
                This function returns an integer (the power domination number of the input graph) when supplied with a NetworkX graph object.
                Optionally, the user may supply this function with a number of compute threads to use in the parallelization step.
                If no number of compute threads are given, then all but one available compute threads are used by default.
            \item
                \verb|PDT_minpds(Input_graph, Number_workers) -> list|
                
                This function returns a list containing vertex labels of vertices that form a minimum power dominating set for the input graph.
                Optionally, the user may supply this function with a number of compute threads to use in the parallelization step.
                If no number of compute threads are given, then all but one available compute threads are used by default.
                The function signature is as follows:

            \item 
                \verb|CheckForPDSOfSize(Input_graph, Contracted, PreferredVertices, CycleEntrances,| \\ \phantom{CheckForPDSOfSize(} \verb|ActiveVertices, Blues, Placement_size, Number_workers) -> list|
                
                This function returns a power dominating set of the given size (if one exists) that is subject to the restrictions discussed in Section 4.
                Power dominating sets returned by this function have no PMUs located on vertices with degree less than 3, no PMUs located on redundant vertices, and PMUs located on all preferred vertices.
                Optionally, the user may supply: a boolean value for if the input graph is already contracted, the list of preferred vertices, the list of paired entrance vertices, the list of active vertices, and the list of vertices colored blue in the restricted power domination problem subject to $\pref{G}$.
                If none of these parameters are supplied, then the PDT will calculate them.
            \item
                \verb|allpdsofsize(Input_graph, Size) -> list|
                
                This generator yields each power dominating set of a given size.
                Note that this function is not parallelized and can be used to find all minimum PDSs and not just ones that satisfy the restricted power domination problem on $G'$ subject to $\pref{G'}$.
            \item 
                \verb|parallel_allpds_of_size(Input_graph, Placement_size, Number_workers)-> list of lists|
                
                This function returns a list of all power dominating sets of a given size.
                This differs from \verb|allpdsofsize| in that each power dominating set is held in memory at a given time.
                Due to memory constraints, this function may overload personal computer's memory capacity for large power dominating sets.
                Note that this function leverages parallel computing methodologies.
        \end{enumerate}
        Documentation for the other functions contained within the PDT can be located on GitHub as well as examples that act as unit tests for each function.
        
        Moreover, we include functions that afford some preliminary investigation into power domination variations within the PDT, including: failed power domination, $k$-fault-tolerant power domination, $k$-PMU-defect-robust power domination, and fragile power domination.
        The PDT can also be used to investigate the token jumping reconfiguration graph for power domination, or the token addition and removal reconfiguration graph for power domination.
        
    \subsection{Acknowledgements}
        The authors thank Dr. Mary Flagg of the University of St. Thomas in Houston for her insights on forts and power domination.
        This project was sponsored, in part, by the Air Force Research Laboratory via the Autonomy Technology Research Center and Wright State University.
        This research was also supported by Air Force Office of Scientific Research award 23RYCOR004 and is Distribution A under the reference number AFRL-2024-1739.

    \subsection{Reproducibility and Code Availability Statement}
        All tests were conducted on a Ryzen 9 5950X with 128 gigabytes of system RAM running Ubuntu 22.04.3.
        Total system RAM usage never exceeded 8 gigabytes at any time during the tests.
        The Power Domination Toolbox is available on GitHub at \url{https://github.com/JibJibFlutterhousen/PowerDominationToolbox}

    \bibliographystyle{plainurl}
    \bibliography{pdt}

\begin{thebibliography}{10}

\bibitem{Anderson2023}
Sarah Anderson, Karen Collins, Daniela Ferrero, Leslie Hogben, Carolyn Mayer,
  Ann Trenk, and Shanise Walker.
\newblock Product throttling for power domination.
\newblock {\em The Australasian Journal of Combinatorics}, 85(3):248--272,
  2023.
\newblock URL: \url{https://ajc.maths.uq.edu.au/pdf/85/ajc_v85_p248.pdf}.

\bibitem{Benson2018}
Katherine Benson, Daniela Ferrero, Mary Flagg, Veronika Furst, Leslie Hogben,
  Violeta Vasilevska, , and Brian Wissman.
\newblock Zero forcing and power domination for graph products.
\newblock {\em Australasian Journal of Combinatorics}, 7(2):221--235, 2018.
\newblock URL: \url{https://ajc.maths.uq.edu.au/pdf/70/ajc_v70_p221.pdf}.

\bibitem{Benson2018a}
Katherine~F. Benson, Daniela Ferrero, Mary Flagg, Veronika Furst, Leslie
  Hogben, and Violeta Vasilevska.
\newblock Nordhaus–gaddum problems for power domination.
\newblock {\em Discrete Applied Mathematics}, 251:103--113, December 2018.
\newblock \href {https://doi.org/10.1016/j.dam.2018.06.004}
  {\path{doi:10.1016/j.dam.2018.06.004}}.

\bibitem{Bjorkman2022}
Beth Bjorkman, Chassidy Bozeman, Daniela Ferrero, Mary Flagg, Cheryl Grood,
  Leslie Hogben, Bonnie Jacob, and Carolyn Reinhart.
\newblock Power domination reconfiguration, 2022.
\newblock \href {https://doi.org/10.48550/ARXIV.2201.01798}
  {\path{doi:10.48550/ARXIV.2201.01798}}.

\bibitem{Bozeman2019}
Chassidy Bozeman, Boris Brimkov, Craig Erickson, Daniela Ferrero, Mary Flagg,
  and Leslie Hogben.
\newblock Restricted power domination and zero forcing problems.
\newblock {\em Journal of Combinatorial Optimization}, 37(3):935--956, July
  2018.
\newblock \href {https://doi.org/10.1007/s10878-018-0330-6}
  {\path{doi:10.1007/s10878-018-0330-6}}.

\bibitem{Erdoes1959}
P.~Erdős and A.~Rényi.
\newblock On random graphs. i.
\newblock {\em Publicationes Mathematicae Debrecen}, 6(3–4):290--297, July
  1949.
\newblock \href {https://doi.org/10.5486/pmd.1959.6.3-4.12}
  {\path{doi:10.5486/pmd.1959.6.3-4.12}}.

\bibitem{Force2004}
Joseph~H Eto.
\newblock Final report on the august 14, 2003 blackout in the united states and
  canada: causes and recommendations, 2004.

\bibitem{Fast2018}
Caleb~C. Fast and Illya~V. Hicks.
\newblock Effects of vertex degrees on the zero-forcing number and propagation
  time of a graph.
\newblock {\em Discrete Applied Mathematics}, 250:215--226, December 2018.
\newblock \href {https://doi.org/10.1016/j.dam.2018.05.002}
  {\path{doi:10.1016/j.dam.2018.05.002}}.

\bibitem{Haynes2002}
Teresa~W. Haynes, Sandra~M. Hedetniemi, Stephen~T. Hedetniemi, and Michael~A.
  Henning.
\newblock Domination in graphs applied to electric power networks.
\newblock {\em SIAM Journal on Discrete Mathematics}, 15(4):519--529, January
  2002.
\newblock \href {https://doi.org/10.1137/s0895480100375831}
  {\path{doi:10.1137/s0895480100375831}}.

\bibitem{hhhh02}
Teresa~W. Haynes, Sandra~M. Hedetniemi, Stephen~T. Hedetniemi, and Michael~A.
  Henning.
\newblock Domination in graphs applied to electric power networks.
\newblock {\em SIAM Journal on Discrete Mathematics}, 15(4):519--529, 2002.
\newblock \href {https://doi.org/10.1137/S0895480100375831}
  {\path{doi:10.1137/S0895480100375831}}.

\bibitem{Hogben2022}
Leslie Hogben, Jephian C.-H. Lin, and Bryan~L. Shader.
\newblock {\em Inverse problems and zero forcing for graphs}.
\newblock Number v. 270 in Mathematical Surveys and Monographs. American
  Mathematical Society, Providence, Rhode Island, 2022.
\newblock Includes bibliographical references and index. - Electronic
  reproduction;Providence, Rhode Island;American Mathematical Society;2022. -
  Description based on print version record.

\bibitem{Baba2021}
Nikolaos~M. Manousakis, George~N. Korres, and Pavlos~S. Georgilakis.
\newblock Taxonomy of pmu placement methodologies.
\newblock {\em IEEE Transactions on Power Systems}, 27(2):1070--1077, May 2012.
\newblock \href {https://doi.org/10.1109/tpwrs.2011.2179816}
  {\path{doi:10.1109/tpwrs.2011.2179816}}.

\bibitem{Graph6}
Brendan McKay.
\newblock URL: \url{http://users.cecs.anu.edu.au/~bdm/data/formats.html}.

\bibitem{IEEE118}
University of~Washington.
\newblock URL: \url{https://labs.ece.uw.edu/pstca/pf118/ieee118cdf.txt}.

\bibitem{Smith2020}
Logan~A. Smith and Illya~V. Hicks.
\newblock Optimal sensor placement in power grids: Power domination, set
  covering, and the neighborhoods of zero forcing forts, 2020.
\newblock \href {https://doi.org/10.48550/ARXIV.2006.03460}
  {\path{doi:10.48550/ARXIV.2006.03460}}.

\bibitem{Thurner2018}
Leon Thurner, Alexander Scheidler, Florian Schafer, Jan-Hendrik Menke, Julian
  Dollichon, Friederike Meier, Steffen Meinecke, and Martin Braun.
\newblock Pandapower—an open-source python tool for convenient modeling,
  analysis, and optimization of electric power systems.
\newblock {\em IEEE Transactions on Power Systems}, 33(6):6510--6521, November
  2018.
\newblock \href {https://doi.org/10.1109/tpwrs.2018.2829021}
  {\path{doi:10.1109/tpwrs.2018.2829021}}.

\end{thebibliography}

    \newpage
    \appendix
    \appendixpage
    \section{Table of IEEE Graph Parameters}
    \begin{table}[H]
        \resizebox{\textwidth}{!}{%
            \begin{tabular}{|l|l|l|l|l|l|l|l|}
                    \hline
                    $|V(G)|$ & \begin{tabular}[c]{@{}l@{}}Number of vertices\\ with $\deg(v) \geq 3$\end{tabular} & $|V(G')|$ & \begin{tabular}[c]{@{}l@{}}Average\\ contraction time\end{tabular} & $|\pref{G'}|$ & \begin{tabular}[c]{@{}l@{}}Average time\\ to calculate\\ $\pref{G'}$\end{tabular} & $|\jeff{G'}|$ & \begin{tabular}[c]{@{}l@{}}Average time\\ to calculate $\jeff{G'}$\end{tabular} \\ \hline
                    4 & 0 & 1 & $4.181\times10^{-5}$ & 0 & $4.321\times10^{-4}$ & 0 & $6.275\times10^{-5}$ \\ \hline
                    5 & 2 & 4 & $8.187\times10^{-5}$ & 0 & $9.654\times10^{-4}$ & 0 & $1.853\times10^{-4}$ \\ \hline
                    6 & 6 & 6 & $5.170\times10^{-5}$ & 0 & $6.929\times10^{-4}$ & 0 & $7.557\times10^{-5}$ \\ \hline
                    9 & 3 & 9 & $9.790\times10^{-5}$ & 0 & $2.265\times10^{-4}$ & 0 & $2.263\times10^{-4}$ \\ \hline
                    11 & 3 & 9 & $1.277\times10^{-4}$ & 2 & $2.525\times10^{-4}$ & 0 & $2.198\times10^{-4}$ \\ \hline
                    14 & 7 & 13 & $1.442\times10^{-4}$ & 0 & $2.088\times10^{-4}$ & 0 & $5.116\times10^{-4}$ \\ \hline
                    24 & 14 & 23 & $2.071\times10^{-4}$ & 0 & $3.165\times10^{-4}$ & 1 & $7.595\times10^{-4}$ \\ \hline
                    30 & 12 & 25 & $2.994\times10^{-4}$ & 1 & $6.621\times10^{-4}$ & 0 & $9.014\times10^{-4}$ \\ \hline
                    30 & 12 & 25 & $3.121\times10^{-4}$ & 1 & $6.611\times10^{-4}$ & 0 & $9.298\times10^{-4}$ \\ \hline
                    33 & 3 & 9 & $3.653\times10^{-4}$ & 2 & $4.501\times10^{-4}$ & 0 & $4.096\times10^{-4}$ \\ \hline
                    39 & 18 & 36 & $3.791\times10^{-4}$ & 3 & $1.293\times10^{-3}$ & 0 & $1.293\times10^{-4}$ \\ \hline
                    57 & 24 & 42 & $6.231\times10^{-4}$ & 0 & $8.031\times10^{-4}$ & 1 & $7.121\times10^{-3}$ \\ \hline
                    89 & 50 & 84 & $1.115\times10^{-3}$ & 3 & $2.973\times10^{-3}$ & 1 & $2.738\times10^{-3}$ \\ \hline
                    118 & 55 & 101 & $1.612\times10^{-3}$ & 1 & $3.422\times10^{-3}$ & 1 & $4.443\times10^{-3}$ \\ \hline
                    145 & 102 & 141 & $2.225\times10^{-3}$ & 4 & $6.811\times10^{-3}$ & 0 & $4.071\times10^{-3}$ \\ \hline
                    200 & 73 & 176 & $3.568\times10^{-3}$ & 17 & $1.022\times10^{-2}$ & 0 & $5.235\times10^{-3}$ \\ \hline
                    300 & 155 & 283 & $6.718\times10^{-3}$ & 11 & $3.818\times10^{-2}$ & 2 & $8.870\times10^{-3}$ \\ \hline
                    1,354 & 496 & 1,233 & $1.119\times10^{-1}$ & 141 & $6.228\times10^{-1}$ & 12 & $1.237\times10^{-1}$ \\ \hline
                \end{tabular}%
        }
        \caption{Parameters for various test networks as available through the pandapower Python module.}
        \label{tab:variousieeegraphs}
    \end{table}

\end{document}